\documentclass[12pt, a4paper]{article}

\usepackage{amssymb}  
\usepackage{amsmath} 
\usepackage{amsthm} 
\usepackage[utf8]{inputenc}
\usepackage[T1]{fontenc}
\usepackage[dvipdfm]{hyperref}

\newtheorem{thm}{Theorem}[section]
\newtheorem{cor}[thm]{Corollary}
\newtheorem{lem}[thm]{Lemma}
\newtheorem{prop}[thm]{Proposition} 
\newtheorem{pte}[thm]{Property} 

\theoremstyle{definition}
\newtheorem{defin}[thm]{Definition}
\newtheorem*{rem}{Remark}
\newtheorem*{convention}{Convention}

\numberwithin{equation}{section}

\newcommand{\R}{\mathbb{R}}
\newcommand{\C}{\mathbb{C}}
\newcommand{\N}{\mathbb{N}}

\newcommand{\g}{\mathfrak{g}}
\newcommand{\ida}{\mathfrak{a}}
\newcommand{\s}{\mathfrak{s}}
\newcommand{\n}{\mathfrak{n}}
\newcommand{\ce}{\mathfrak{t}} 

\newcommand{\schw}{\mathcal{S}}
\newcommand{\SG}{\mathcal{\schw}_\sigma(G)}
\newcommand{\SGi}{\mathcal{\schw}_{\sigma_1}(G_1)}
\newcommand{\SGii}{\mathcal{\schw}_{\sigma_2}(G_2)}

\newcommand{\opg}{{\textrm{op}\,\g}}
\newcommand{\opgi}{{\textrm{op}\,\g_i}}
\newcommand{\opn}{{\textrm{op}\,\n}}

\newcommand{\Id}{\operatorname{Id}}
\newcommand{\Ad}{\operatorname{Ad}}
\newcommand{\ad}{\operatorname{ad}}
\newcommand{\vect}{\operatorname{Vect}}
\newcommand{\supp}{\operatorname{supp}}

\newcommand{\indUl}{\ensuremath{\mathbf{1}_{U^l}}}
\newcommand{\indUlo}{\ensuremath{\mathbf{1}_{U^{l_0}}}}\newcommand{\pdtN}{\ensuremath{\cdot_\textup{\tiny CBH}}}

\begin{document}


\title{Schwartz functions, tempered distributions, and Kernel Theorem on solvable Lie groups}

\author{Emilie David-Guillou}

\date{February 10, 2010}

\maketitle


\begin{abstract}
Let $G$ be a solvable Lie group endowed with right Haar measure. We define and study a dense Fréchet $\ast$-subalgebra $\SG$ of $L^1(G)$, consisting of smooth functions 
rapidly-decreasing at infinity on $G$. 
When $G$ is nilpotent, we recover the classical Schwartz algebra $\schw(G)$. We develop  a distribution theory for $\SG$, and we generalize the classical Kernel Theorem of L.~Schwartz to our setting. 
\end{abstract}


\section{Introduction}

The aim of the present paper is to define an analogue of Schwartz functions and tempered distributions on connected simply connected solvable Lie groups, and to generalize some classical results from Schwartz theory to that setting. 

Let $G$ be a connected simply connected solvable Lie group endowed with the right invariant Haar measure. 
We are interested in forming a topological \mbox{$\ast$-algebra} of smooth integrable functions on $G$, without assumption of compactness on the supports. Practically, we look for an intermediate function algebra between $L^1(G)$ and the algebra $C_0^\infty(G)$ of smooth functions with compact support, that would be modelled on the Euclidean Schwartz algebra. In other words, we look for an algebra of smooth functions that decrease, with all their derivatives, rapidly at infinity on $G$. 

If $G$ is nilpotent there is a naturally defined Schwartz space $\schw(G)$ (see \emph{e.g.} 
\cite{howe77} or \cite{corwin81}), which has  these characteristics. In this case, indeed, the exponential map $\exp: \g \to G$ is a diffeomorphism from the Lie algebra $\g$ of $G$ onto $G$. This allows to identify the Lie group $G$ and its Lie algebra $\g$, and the product law on $G=\g$ is then given by the Campbell-Baker-Hausdorff product
\[ X\cdot_{\textrm{\tiny CBH}}Y=X+Y+\frac{1}{2}[X,Y]+\dots \, , \qquad X,Y\in \g ,  \]
which is in fact a polynomial mapping since $G$ is nilpotent. One defines the Schwartz space $\schw(G)$
as the image under composition with the exponential map of the usual Schwartz space $\schw(\g)$ of rapidly decreasing smooth functions on $\g$ (seen as finite dimensional real vector space). The topology on $\schw(G)$ is that for which the composition with the exponential map is a topological isomorphism of $\schw(\g)$ into $\schw(G)$.  
The Schwartz space $\schw(G)$ has the properties that we want by construction: it consists of smooth functions, 
it contains $C^\infty_0(G)$ as dense subspace, and it is a dense \mbox{$\ast$-subalgebra} of $L^1(G)$.

If $G$ is not nilpotent but only solvable, the ``natural way'' to define a Schwartz space on $G$ is less clear. 
Even if we assume that the exponential map $\exp: \g \to G$ is a diffeomorphism (in which case $G$ is said to be an \emph{exponential} solvable Lie group), there is no hope to define $\schw(G)$  as the image of the Euclidean Schwartz space $\schw(\g)$ under composition with the exponential map, because there are functions $f$ in $\schw(\g)$ such that $f\circ \exp^{-1}$ fails to be integrable on $G$. The problem on such a group $G$, is that the volume of the balls grows exponentially with respect to their radius, so the decay imposed on the function $f$ is not rapid enough to imply the integrability of $f\circ \exp^{-1}$ on $G$. 

In this paper, we define a property of rapid decay on a general connected simply connected solvable $G$, that takes into account the geometry of $G$, by asking exponential decay 
in certain directions. With respect to this notion, we form a space
of smooth functions on $G$ that decrease with their derivatives rapidly at infinity. 
We endow this space with a family of seminorms, and 
we obtain a dense Fréchet $\ast$-subalgebra of $L^1(G)$ (Theorem~\ref{exstnce_pds'}), with the $C^\infty$ functions with compact support as dense $\ast$-subalgebra (Theorem~\ref{DSdense}).  
If $G$ is nilpotent, our algebra of rapidly decreasing smooth functions agrees with the Schwartz space $\schw(G)$ (Proposition~\ref{SG=Snilp}). 
The generality of our construction is discussed in \S\ref{arbitrary}. We give explicit descriptions of the algebra $\SG$ for some typical classes of solvable Lie groups $G$ in \S\ref{exemples}. 

Next, we investigate the properties of the dual space of our space of rapidly decreasing smooth functions. We find that every element of the dual space is a finite linear combination of slowly increasing functions and derivatives of such functions  (Theorem~\ref{D-temperee}). 
In this respect,  
we retrieve the classical concept of tempered distributions.

The last aspect considered in the paper in that of  product groups. We show that the space of rapidly decreasing functions on the direct product  $G_1\times G_2$ of two connected simply connected solvable Lie groups, is the completed tensor product of $\SGi$ and $\SGii$ for the projective topology (Theorem~\ref{Stensor}). As corollary, we obtain a variant of Schwartz Kernel Theorem. Recall that in its Euclidean version,  Schwartz Kernel Theorem asserts essentially that every continuous linear map from the space of rapidly decreasing functions $\schw_x(\R^n)$ into the space of tempered distributions  $\schw'_y(\R^m)$, is given by a unique distribution ``in both variables $x$ and $y$'', and that the correspondence between linear forms and distributions is of topological nature. 
We prove a similar result for our classes of rapidly decreasing functions on connected simply connected solvable Lie groups (Corollary~\ref{skt}).

\section{Towards Schwartz algebra -- Preliminaries}
\label{preliminaires}

\subsection{Theorem of existence} \label{algebre_de_S}

Let $G$ be a connected simply connected solvable Lie group endowed with right invariant Haar measure $dg$, and $\g$ be its Lie algebra. Denote by $m$ the dimension of $\g$. 
Fix $\{X_1, \ldots, X_m\}$ a basis 
of $\g$. To each $\alpha=(\alpha_1, \ldots, \alpha_m)\in \N^m$, we associate a left invariant differential operator $X^\alpha$ on $G$, by setting  for $\phi \in C^\infty(G)$  
\[ X^\alpha \phi (g) = \left. \frac{\partial^{\alpha_1}}{\partial t_1^{\alpha_1}} \ldots \frac{\partial^{\alpha_m}}{\partial t_m^{\alpha_m}} \phi\big(g\exp(t_1X_1)\dots\exp(t_mX_m)\big)\right|_{t_1=\ldots=t_m=0} \, .   \]

\begin{defin} \label{def_S}
Let $\varsigma$ be a Borel function with values in $\R^+$; we say that $\varsigma$ is a \emph{weight function} on $G$.  
Define $\schw_\varsigma(G)$ as being the set of $C^\infty$   
functions $\phi$ on $G$ such that 
\[  
\| \phi \|^{\infty}_{k,\alpha}:= \| \varsigma^k X^\alpha\phi \|_{L^\infty(G)} <\infty 
\]
for every $k\in \N$ and $\alpha \in \N^m$.
Clearly $\schw_\varsigma(G)$ 
has a vector space structure. We call $\schw_\varsigma(G)$ the \emph{space of smooth functions decreasing $\varsigma$-rapidly at infinity on $G$}, and we equip  
$\schw_\varsigma(G)$
with the topology of the seminorms $\| \cdot \|^\infty_{k,\alpha}$. 
\end{defin}

\begin{rem}
The definition of the (topological) vector space $\schw_\varsigma(G)$ is independent of the choice of the basis $\{X_1, \ldots, X_m\}$ because
all the left invariant differential operators are finite linear combinations of   $X^\alpha$'s by Poincaré-Birkhoff-Witt theorem. 
\end{rem}

By construction  $\schw_\varsigma(G)$ is locally convex,  
but an adequate choice of $\varsigma$ makes its structure much richer: 

\begin{thm} \label{exstnce_pds}
Let $G$ and $dg$ be as above. 
There exists a weight function $\varsigma$ on $G$, 
such that $\schw_\varsigma(G)$ is a dense Fréchet $\ast$-subalgebra of the convolution algebra $L^1(G, dg)$.
\end{thm}

In order to prove Theorem~\ref{exstnce_pds},  
we construct an explicit weight function $\varsigma$ with the required property.  
For this purpose, it is convenient to consider some particular realization for $G$, in which the weight function $\varsigma$ has a simple expression. 
The construction  is well known by experts (see \emph{e.g.} \cite{alexopoulos92} or \cite{leptin-ludwig94}).    
For sake of completeness,  we give the main lines 
below; the notations are those from \cite{inoue-ludwig07}. 

\subsection{Special realization of the Lie group $G$} \label{decomp_G}

Let $\n$ be the nilradical of $\g$. 
A linear form $\lambda$ on the complexification $\g_\C$ of the Lie algebra $\g$, is called a \emph{root} if there exists an ideal $\ida$ of $\g_\C$ and some non-zero $w\in \g_\C \setminus \ida $, 
such that $\ad(v)w-\lambda(v)w\in \ida$ for all $v\in \g$. We denote by $\mathcal{R}$ the set of roots of $\g$. 
Take $X\in \g$ \emph{in general position with respect to 
$\mathcal{R}$} (which means that $\lambda(X)\neq\lambda'(X)$ for any two distinct roots $\lambda$, $\lambda'$). Then 
$\{\lambda(X); \, \lambda\in \mathcal{R}\}$
is the set of all the distinct eigenvalues of the linear operator $\ad(X)$ on $\g_ \C$, and thus we have a Jordan decomposition of the space $\g_\C$ 
\[\g_\C =
\bigoplus_{\lambda\in \mathcal{R}} 
\g_{\C, \lambda}  \]
along the generalized eigenspaces of $\ad(X)$
\[\g_{\C,\lambda} =\big\{ Y\in \g_\C;\, 
(\ad(X)  -\lambda(X)
  \Id_{\g_\C})^{\dim \g}(Y)=0 \big \}. \]  
By standard arguments, we have
\[[\g_{\C,\lambda},\g_{\C,\lambda'} ]\subset 
\g_{\C,\lambda+\lambda'} 
\qquad \lambda, \lambda'\in \mathcal{R},\]  
which implies in particular that $\g_0:=\g_{\C, 0}\cap\g$ is a nilpotent subalgebra of $\g$. Furthermore 
\[\bigoplus_{\lambda\in \mathcal{R}\setminus\{0\}} \g_{\C, \lambda}\subset
[\g_\C,\g_\C]\subset \n_\C, \] 
where $\n_\C$ denotes the complexification of the nilradical $\n$. 
It follows that 
\[\g=\g_0+\n.\] 
We choose a subspace $\ce$ in $\g_0$ such that \[ \g=\ce \oplus\n.\]
We will now 
define a Lie group structure on
$\boldsymbol{G}:=\ce\times\n$.
Since $\ce\subset \g_0$ and $\g_0$ nilpotent, the Campbell-Baker-Hausdorff
product of two elements $T$ and $T'$ of $\ce$ is given by a finite expression
\[T\cdot_\textrm{\tiny CBH}T'=T+T'+\frac{1}{2}
[T,T']+\dots=T+T'+P(T,T')\]
where $P:\ce\times\ce \to \g_0\cap\n$  is a polynomial mapping. We
endow $\boldsymbol{G}$ with the product law in coordinates $(t,x)\in\ce\times\n$ 
\begin{equation} \label{pdt_sur_Q}
(t,x)\cdot(t',x')=(t+t',P(t,t') \, \cdot_\textrm{\tiny CBH} x \,
\cdot_\textrm{\tiny CBH} e^{\ad t}x').
\end{equation}
It is easy to check that, equipped with this product,
$\boldsymbol{G}$ is a Lie group which admits as Lie algebra $\ce\times\n$ with the Lie bracket 
\begin{equation} \label{s_bracket} [(T,X),(T',X')]=(0,[T,T']+[X,X']+[T,X']+[T',X]).  \end{equation}
In particular, the Lie algebra of $\boldsymbol{G}$ is isomorphic to $\g$. Since $\boldsymbol{G}$ is connected and simply connected by construction, it means that as Lie group $\boldsymbol{G}=G$. 

\begin{convention} In the rest of paper, $\n$, $\ce$, and $P:\ce\times\ce \to \g_0\cap\n$ will be as above,
and we will use 
the realization $(\boldsymbol{G},\, \cdot)$ for the Lie group 
$G$. So from now on,  $G=\ce\times\n$  with the product law~(\ref{pdt_sur_Q}). 
\end{convention}

\begin{rem}
The Lie algebra of the analytic subgroup $N=\{0\}\times \n$ of $G$ is isomorphic to $\n$. Since $N$ is simply connected,  $N$ is by definition the nilradical of $G$. In what follows, we will use the notation $N$ equally for the subgroup $\{0\}\times \n$ of $G$ and for the connected simply connected nilpotent Lie group $(\n, \, \cdot_\textup{\tiny CBH} )$.
\end{rem}

\subsection{More definitions and notations} 
\label{notations}

Let $d$, resp. $k=m-d$ be the dimension of $\n$, resp. $\ce$, as real vector space. Since $G$ is solvable, $d \ge 1$ (while $k$ can be zero). 
We denote by $dn$, resp.~$dt$, the Lebesgue measure on $\n$, resp.~$\ce$. 
Then $dg=dt\,dn$ is a right invariant Haar measure on $G$, and $d^lg= e^{-\ad t} dt \,dn$ the corresponding left invariant one. We write $\delta(g)=\delta(t,n)= e^{ \ad t}$ the modular function on $G$.  

Denote by $| \cdot |_{\R^d}$, $| \cdot |_{\R^k}$, $| \cdot |_{\R^m}$ the Euclidean norms on $\n$, $\ce$, $\g = \ce \times \n$ respectively, and by $\| \cdot \|_{\textrm{op}\,\n}$, $\| \cdot \|_{\textrm{op}\,\g}$ the related operator norms on $\n$, $\g$.

Fix $U$ some symmetric compact neighbourhood of identity $e$ in $G$. We define the \emph{length} of an element $g$ of $G$ 
by
\[
 |g|_G = \inf\{j\in \N; \ g\in U^j= U \cdots U
\} \]
where by convention $U^0=\{e\}$. It is easy to check that two different symmetric compact neighbourhoods $U$ and $U'$ of identity define equivalent lengths in the sense that, for some $C>0$, 
\begin{equation} \label{equiv-lgueurs} 
C^{-1}|g|_{U'} \le |g|_U \le C |g|_{U'}, \qquad g\in G.\end{equation}
Note also that the length is a subadditive function on $G$ 
\begin{equation} \label{l_ssadtve}
|g\cdot h|_G \le |g|_G + |h|_G, \qquad g,h\in G, 
\end{equation} 
and since $U$ is symmetric, that it satisfies
\begin{equation}  \label{l_symtric}
 |g^{-1}|_G=|g|_G, \qquad g\in G. 
\end{equation}

Since $G$ is connected, the nilradical $N$ of $G$ is a closed analytic subgroup of $G$ and $V=U\cap N$ is a symmetric compact neighbourhood of identity in $N$. In addition to  the length $|n|_G$ of elements $n$ of $N$ seen as elements of $G$, we endow $N$ with an intrinsic length
\[ |n|_N = \inf\{j\in \N; \ n\in V^j= V \cdots V\}\]
The two above lengths on $N$ are related in the following way
\begin{equation} \label{comparaison_QN} 
|n|_G \le |n|_N \le C \exp (C |n|_G), \qquad n\in N. 
\end{equation}
The left hand inequality is clear from the definitions of $| \cdot|_N$ and $|\cdot|_G$. The right hand inequality asserts that the nilradical of a connected simply connected (solvable) Lie group has at most exponential distortion (see \cite{varopoulos94b}).  

For a function $\phi$ on $G$,  we use the notation
$\check{\phi}$ for 
\[ \check{\phi}(g):= \phi(g^{-1}), \qquad g\in G. 
\]

\section{Rapidly decreasing functions on solvable Lie groups}
\label{S_algebra} 

\subsection{Reformulation of Theorem~\ref{exstnce_pds} with explicit 
weight function}

Let $\sigma$ be the weight function defined at  $g=(t,n)\in G$ by 
\begin{equation} \label{def_sigma}
\sigma(g):=   \max\big( \|\Ad (g)\|_\opg, \,  \|\Ad(g^{-1})\|_\opg \big) \cdot (1+|g|_G+|n|_N)
, \end{equation} 
where $\Ad$ denotes the adjoint representation of $G$ on its Lie algebra  $\g$. 
We reformulate Theorem~\ref{exstnce_pds} as follows: 
\begin{thm} \label{exstnce_pds'} 
Let $G$, $dg$ be as in \S\ref{preliminaires}, and let $\sigma$ given by~(\ref{def_sigma}).  Then the space $\schw_\sigma(G)$ of smooth functions decreasing $\sigma$-rapidly at infinity on $G$,  is a dense Fréchet $\ast$-subalgebra of $L^1(G, dg)$. 
\end{thm} 

\begin{rem}
The space $\SG$ is independent of the choice of the compact neighbourhood $U$, hence of the definitions of the lengths $|\cdot|_G$ and $|\cdot|_N$. Changing $U$ into $U'$ changes  the weight function $\sigma$ into 
$\sigma'$ with $c \sigma' \le \sigma \le C \sigma'$ by~(\ref{equiv-lgueurs}), and thus
$\SG=\schw_{\sigma'}(G)$. 
\end{rem}

\subsection{Properties of $\sigma$}

We give here some properties of the weight function $\sigma$, that will be used to prove {Theorem~\ref{exstnce_pds'}} in \S\ref{pf_1'}. 

\begin{pte} \label{compense_volume} 
There exists $p \in \N$ such that 
\[ \int_G \frac{1}{\sigma^p(g)} dg < \infty \,.\]
We say that $\sigma$ \emph{compensates the growth of the volume of $G$}.
\end{pte}

\begin{proof} 
Since the Haar measure $dg$ on $G$ splits into $dt\,dn$, we can treat separately the volume growth on $\ce$ and the volume growth on $\n$. 

In the $\ce$ coordinate, the growth of the volume is compensated by $|g|_G$, thanks to the estimate
\begin{equation} \label{controle_g} 
|t|_{\R^k} \le C |g|_G   , \qquad g=(t,n) \in G.  
\end{equation}
Let us prove~(\ref{controle_g}). Assume that $|g|_G=j$. Then by definition, there are $j$ elements $(t_1,n_1), \ldots , 
(t_j, n_j)$ in the compact subset $U$, 
such that 
$g = (t_1,n_1) \cdots (t_j,n_j)$.
Now the product law on $G=(\ce\times \n)$ 
is just a sum in the $\ce$ component, therefore
\[ |t|_{\R^k} =|t_1+ \ldots + t_j|_{\R^k} 
\le j  \sup\{|a|_{\R^k}; \, (a,b)\in U\} = C |g|_G, \qquad g=(t,n). \]

In the $\n$ coordinate, the growth of the volume is compensated by $|n|_N$. Indeed, recall that $N=(\n, \pdtN)$. Then it is easy to see that 
$|n|_{\R^d}$ is dominated by some power of the length $|n|_N$ (see \cite{dixmier60}). 
In fact, the reverse estimate is also valid, and 
so there are  $q\in \N^\ast$ and $C>1$ such that   
\begin{equation} \label{controle_u}
C^{-1}\, |n|_{\R^d}^{1/q}  \, \le \,  |n|_N \, \le \,  C \,  \big(|n|_{\R^d}+1\big), \qquad n\in N   
\end{equation}
(see \emph{e.g.} \cite{ludwig-molitorbraun95}).

Putting~(\ref{controle_g}) and~(\ref{controle_u}) together shows
\[ \sigma(g) \ge 1+|g|_G+|n|_N \ge C (1+|t|_{\R^k}+|n|_{\R^d}^{1/q}) \ge C (1+|t|_{\R^k}+|n|_{\R^d})^{1/q}
, \]
with $C>0$ independent of $g=(t,n)\in G$. 
It follows that, if $p\ge \frac{k+d+2}{q}$,   
\begin{align*}
\int_G \frac{1}{\sigma^{p}(g)} \,dg &\le 
C \int_{\R^k} \int_{\R^d}  \frac{1}{(1+|t|_{\R^k}+|n|_{\R^d})^{k+d+2}} \,dn dt \\ 
&\le C \int_{\R^{k+d}}   \frac{1}{(1+|(t,n)|_{\R^{k+d}})^{k+d+2}} \,dtdn \quad  <  \quad  \infty. \qedhere
\end{align*}
\end{proof}

\begin{pte} \label{controle_module}
The modular function of $G$  
is dominated by a power of the weight function $\sigma$. For all $g\in G$, we have 
\[ \delta(g) \le  \sigma^m(g) .\]
\end{pte}

\begin{proof} 
Right and left Haar measures $dg$ and $d^lg$ are obtained by translating a nonvanishing differential form $\omega \in \bigwedge^{m} T^\ast_eG$.  
The action of $G$ on $\bigwedge^{m}T^\ast_eG $
induces at a point $g\in G$ a corrective factor between the two measures, the modular function of the group, which is 
$\delta(g) = | \det \Ad(g)|.$
Since the operator norm dominates the module of every eigenvalue, we have $| \det \Ad(g) |\le \| \Ad(g) \|_\opg^{m}$. It follows that 
\[\delta(g) \le \| \Ad(g) \|_\opg^{m} \le \sigma^m(g), \qquad g\in G. \qedhere \] \end{proof}

\begin{pte} \label{equiv_inverse} 
Each one of the  
functions $\sigma$ and $\check{\sigma}$ is
dominated by a power of the other one, up to some constant factor.  
Summarized in one equation: there are $C>1$ and $r\in \N^\ast$ such that \[ C^{-1} \sigma^{1/r}(g) \le \check{\sigma}(g) \le C \sigma^{r}(g), \qquad g\in G. \] 
\end{pte}

\begin{proof}
The inequality on the left hand side follows from the one on the right hand side by change of variable $g \leftrightarrow g^{-1}$, so it suffices to show that 
$\check{\sigma}$ is dominated by a power of  $\sigma$.
Let $g=(t,n)\in G$. By definition 
\begin{align*}
\check{\sigma}(g) &= \sigma(g^{-1})=\sigma(-t,-e^{\ad t} n) \\ 
&=  \max\big( \|\Ad (g^{-1})\|_\opg, \,  \|\Ad(g)\|_\opg \big) \cdot \big(1+|g^{-1}|_G+|-e^{\ad t} n|_N\big).  
\end{align*}
In the expression above, \begin{itemize}
\item the factor on the left is unchanged by the transformation $g\mapsto g^{-1}$; 
\item $|g^{-1}|_G=|g|_G$ by~(\ref{l_symtric})
\end{itemize}
So to dominate $\check{\sigma}$ by some power of $\sigma$, we need only to bound 
$|-e^{\ad t} n|_N$ adequately. 

In $(N, \, \pdtN)$, the inverse of an element $n$ is $-n$, so $|-e^{\ad t} n|_N= |e^{\ad t} n|_N$. We apply~(\ref{controle_u}), and find that for all  $(t,n)\in \ce\times\n$
\begin{equation} \label{e^adt}
|e^{\ad t} n|_N  \le C  \big(|e^{\ad t} n|_{\R^d}+1\big) \le  C  \big(\|e^{\ad t}\|_\opn \cdot | n|_{\R^d}+1\big).
\end{equation}

From the construction in \S\ref{decomp_G} follows that the operator $e^{\ad t}:\n\to\n$ 
is 
the restriction to $\n$ of the adjoint operator $\Ad (t,0)$ with  $(t,0)\in \ce\times\n$ (the restriction is well defined because, as nilradical, $\n$ is an ideal of the Lie algebra $\g$).
Writing $(t,0)=(0,-n)\cdot (t,n)$, we see that  $\Ad(t,0)=\Ad(0,-n)\circ \Ad(t,n)$ as operators on $\g$. By restriction to the subspace $\n$, we obtain that $e^{\ad t}=\Ad(0,-n)\vert_\n\circ \Ad(t,n)\vert_\n$.

By construction again, the restriction to $\n$ of the operator $\Ad(0,n)$ coincides with the adjoint operator $\Ad(n)$ on the Lie group $N$. Since $N =(\n,\, \pdtN)$, the exponential mapping from $\n$ onto $N$ is the identity map. So $\Ad(n)=e^{\ad(n)}=P(\ad(n))$ for some real polynomial $P$, because $\ad(n)$ is a nilpotent operator on $\n$. By linearity of $\ad$, there are $C>0$ and $q'\in \N$ such that 
\begin{equation} \label{domination_Ad}
\| \Ad(n)\|_{\opn}  \le \| P(\ad(n)) \|_{\opn} \le C(|n|_{\R^d}^{q'}+1), \qquad n\in N.  \end{equation} It follows that
\begin{align} \label{op_e^adt}
\|e^{\ad t} \|_\opn &\le \big\|\Ad(0,-n)\vert_\n \big\|_\opn \cdot \big\| \Ad (t,n)\vert_\n \big\|_\opn  \nonumber 
\\  
&\le \|\Ad(-n) \|_\opn \cdot \| \Ad (t,n) \|_\opg 
\nonumber \\
&\le C (|n|_{\R^d}^{q'}+1) \cdot \| \Ad(g) \|_\opg, 
\end{align}
with $C>0$ independent of $g=(t,n)\in G$. 

By injecting this estimate in~(\ref{e^adt}), then  using~(\ref{controle_u}), we obtain that 
\begin{align*} 
|e^{\ad t} n |_N &\le C \big( (|n|_{\R^d}^{q'}+1) \cdot \| \Ad(g) \|_\opg \cdot |n|_{\R^d} +1 \big) \\ 
&\le C \big( (|n|_{N}^{q(q'+1)}+1) \cdot \| \Ad(g) \|_\opg  +1 \big), \qquad g=(t,n)\in G.  
\end{align*}
It implies that 
\begin{align*}
\check{\sigma}(g) &\le C \max\big( \|\Ad (g)\|_\opg, \,  \|\Ad(g^{-1})\|_\opg \big)^2 \cdot \big(1+|g|_G+|n|^{q(q'+1)}_N\big) \\ &\le C \sigma(g)^{\max(2,q(q'+1))}, 
\end{align*}
for all  $g=(t,n)\in G$. 
This proves  Property~\ref{equiv_inverse} with $r={\max(2,q(q'+1))}$.
\end{proof}

\begin{pte} \label{sous-polynomial}
There exist $C>0$ and $s\in \N$ such that 
\[\sigma(g\cdot g')\le C \sigma^s(g) \cdot
 \sigma^s(g'), \qquad g,g'\in G. \]
We say that $\sigma$ is \emph{sub-polynomial}. 
\end{pte}

\begin{proof} Let $g=(t,n), \, g'=(t',n') \in G$. By definition, 
\begin{align*} 
\sigma(g\cdot g')  &=   \max\big(\| \Ad(g\cdot g')\|_\opg,  \| \Ad( (g\cdot g')^{-1})\|_\opg \big)  \nonumber \\
&\qquad \times (1+|g\cdot g'|_G+|n+e^{\ad t} n'|_N) . \end{align*}
In this expression, some terms 
are clearly compatible with sub-polynomiality
\begin{itemize}
\item 
$\max\big( \|\Ad(g)\|_\opg, \,  \|\Ad(g^{-1})\|_\opg\big)$, because as function of $g\in G$ it is submulticative 
(which is clear since $\Ad$ is an homomorphism of $G$ into the group of automorphisms on $\g$);
\item $1+|g\cdot g'|_G$, because,  by ~(\ref{l_ssadtve}), $ 1 +|g\cdot g'|_G \le (1+|g|_G)  (1+|g'|_G)$. 
\end{itemize}
So we just need to  control the term $|n+e^{\ad t} n'|_N$. 
Using the estimates~(\ref{controle_u}) and~(\ref{op_e^adt}), we have that 
\begin{align*}
|n+e^{\ad t} n'|_N &\le C ( |n+e^{\ad t} n'|_{\R^d}+1) \\ 
&\le C ( |n|_{\R^d}+\| e^{\ad t} \|_\opn \cdot |n'|_{\R^d}+1) \\
&\le C ( |n|_{\R^d}+ (|n|^{q'}_{\R^d} +1) \cdot \| \Ad(g) \|_\opn \cdot |n'|_{\R^d}+1) \\
&\le C ( |n|_N^q + (|n|^{qq'}_N +1) \cdot \| \Ad(g) \|_\opn \cdot |n'|^q_{N}+1), 
\end{align*}
for all $g=(n,t)$ and $g'=(n',t')$ in $G$. 

Back to $\sigma(g\cdot g')$, this implies \begin{align*}
\sigma(g\cdot g') &\le 
\max \big(\| \Ad(g)\|_\opg,  \| \Ad( g^{-1})\|_\opg \big)^2 (1+|g|_G +|n|_N)^{q\max(1,q')} \\
&\quad \times \max\big(\| \Ad(g')\|_\opg,  \| \Ad( g'^{-1})\|_\opg \big)  \cdot (1+|g'|_G +|n'|_N)^q \\
&\le C \sigma(g)^{\max(2, q(q'+1))} \cdot \sigma(g')^q, 
\end{align*}
with $C>0$ and $q,q'\in \N$ independent of $g$ and $g'$. 
\end{proof}

\subsection{Proof of Theorem~\ref{exstnce_pds'}} \label{pf_1'}

\subsubsection*{Step one: $\schw_\sigma(G)$ is a dense subset of $L^1(G,dg)$.}

The inclusion of $\schw_\sigma(G)$ in  $L^1(G,dg)$ follows easily from the fact that $\sigma$ compensate the growth of the volume of the group $G$ (Property~\ref{compense_volume}). Indeed if $\phi\in \schw_\sigma(G)$, 
\[ \int_G |\phi(g)| dg \le  \int_G \frac{1}{\sigma^{p}(g)} dg \cdot \| \sigma^p \phi\|_{L^\infty(G)} = C  \| \phi \|^\infty_{p,0} < \infty,   \]
so $\phi\in L^1(G, dg)$. 

The density of $\schw_\sigma(G)$ is clear since the set $C^\infty_0(G)$ of differentiable compactly supported functions on $G$ is dense in $L^1(G, dg)$, and it is obvious that $C^\infty_0(G) \subset \schw_\sigma(G)$. 
 
\subsubsection*{{Step two: $\schw_\sigma(G)$ is a Fréchet space.}}

The family of seminorms $\| \cdot \|^{\infty}_{k,\alpha}$, $k\in \N$, $\alpha \in \N^m$, is countable and separating  ($\|\phi\|^\infty_{0,0}=\|\phi\|_{L^\infty(G)}=0 \Rightarrow
\phi = 0$) on the vector space $\schw_\sigma(G)$,  
so it defines a locally convex topology on $\schw_\sigma(G)$. To prove that  $\schw_\sigma(G)$ is a Fréchet space, we just have to show that it is complete. The proof is modelled on the Euclidean case.  

Let $\{\phi_n\}_{n\in \N}\subset \schw_\sigma(G)$ be a Cauchy sequence for the seminorms $\| \cdot \|_{k,\alpha}^{\infty}$. 
For every $k\in \N$ and $\alpha\in \N^m$, the functions $\sigma^k X^\alpha \phi_n$ converge uniformly to a bounded function $\phi_{k,\alpha}$. 
If we show that 
\begin{equation} \label{phi_k,alpha}
\phi_{k,\alpha}=\sigma^k X^\alpha\phi_{0,0}, \qquad k\in\N, \, \alpha \in \N^m,  
\end{equation}
it will prove that $\phi_{0,0}$ belong to the space $\schw_\sigma(G)$ and that $\phi_n$ converges to $\phi_{0,0}$ in $\schw_\sigma(G)$, and thus that $\schw_\sigma(G)$ is complete.  

Let us prove~(\ref{phi_k,alpha}).  
For $k=0$ and $\alpha$ of length one, say $\alpha=\alpha_i$ with all coordinates equal to zero but the $i^\textrm{th}$ equal to one, we have for all $t\in \R$ 
\[ \phi_n(g\exp(t X_i)) = \phi_n(g) +  \int_0^t X_i \phi_n(g\exp(s X_i)) ds\,.  \] 
Letting $n$ go to infinity in the above equality gives 
\[ \phi_{0,0}(g\exp(t X_i)) = \phi_{0,0}(g) +  \int_0^t  \phi_{0,\alpha_i}(g\exp(s X_i)) ds\, . \]
We differentiate with respect to $t$ at 0. It shows that $\phi_{0,0}$ is continuously differentiable in the direction $X_i$ with
\[ X_i\phi_{0,0}(g) = \phi_{0,\alpha_i}(g), \qquad g\in G.  \] 
Repeating the argument shows that $\phi_{0,0} \in C^\infty(G)$ with $X^\alpha \phi_{0,0}= \phi_{0,\alpha} \ \forall \alpha\in \N^m$. 
This implies that for all $k\in \N$ and $\alpha\in \N^m$, $\sigma^k X^\alpha \phi_n$ converges pointwise to $\sigma^k X^\alpha \phi_{0, 0}$. 
But $\sigma^k X^\alpha \phi_n$ converges uniformly to $ \phi_{k, \alpha}$ by hypothesis, so $\phi_{k, \alpha}= \sigma^k X^\alpha \phi_{0, 0}$. 
Since $\phi_{k, \alpha}\in L^\infty(G)$, this proves that so $\phi_{0,0}\in \schw_\sigma(G)$ and that 
$\phi_n$ converges to $\phi_{0,0}$ in $\schw_\sigma(G)$. Thus $\schw_\sigma(G)$ is complete, and therefore Fréchet. 

\subsubsection*{{Step three: Convolution is continuous from $\schw_\sigma(G) \times \schw_\sigma(G)$ to $ \schw_\sigma(G)$.}}

The convolution of two measurable functions $\phi, \psi$ on $G$,  is defined by 
\[\phi \ast  \psi (g) = \int_G \phi(gh^{-1})\psi(h) \, dh = \int_G \phi(h) \psi(h^{-1}g) \, d^lh, \qquad g\in G, \] provided that the integrals converge. It is easy to check that since the differential operators $X^\alpha$ are left invariant, they act on the convolution in the following way:
\[ X^\alpha (\phi \ast \psi) = \phi  \ast X^\alpha \psi. \]

Let $\phi,\, \psi \in \schw_\sigma(G)$, $k\in \N$ and $\alpha \in \N^m$. By Property~\ref{sous-polynomial} and Property~\ref{compense_volume}, 
\begin{align*}
\|\phi\ast\psi \|^\infty_{k,\alpha} &= \sup_{g\in G} \left| \sigma^k(g) \int_G \phi(gh^{-1}) X^\alpha \psi (h) \,dh \right| \\
&\le \sup_{g\in G}  \left| C \int_G \sigma^{ks}(gh^{-1}) \phi(gh^{-1}) \sigma^{ks}(h) X^\alpha \psi (h) \,dh \right|
\\
&\le C \| \phi \|^\infty_{ks,0} \int_G \frac{ \sigma^{ks+p}(h) X^\alpha \psi (h)}{\sigma^{p}(h)} \,dh \\
&\le C \| \phi \|^\infty_{ks,0} \| \psi\|^\infty_{ks+p,\alpha} ,  
\end{align*}
hence the continuity. 

\subsubsection*{Step four: $L^1(G,dg)$-involution is continuous on  $\schw_\sigma(G)$} 

The involution on $L^1(G,dg)$ is defined by $\phi^\ast(g)= \overline{\phi(g^{-1})} \delta(g^{-1})$.

For $\psi\in C^\infty(G)$ and $\alpha=\alpha_j$ (with the notation of \emph{Step Two}), 
\begin{align*} 
X^\alpha \check{\psi} (g) 
&= \left. \frac{d}{dt} \check{\psi}\big(g\exp(tX_j)\,\big)\right|_{t=0}
= \left. \frac{d}{dt} \psi\big(\exp(-tX_j)g^{-1}\big)\right|_{t=0} \\ 
& = \left. \frac{d}{dt}
\psi\big(g^{-1}\exp(-t\Ad(g)X_j)\big)\,\right|_{t=0} 
= -(\Ad(g)X_j \psi)(g^{-1}) \\
&= - \sum_{1\le i \le m} (\Ad(g)_{ij} X_i\psi) (g^{-1}),  \qquad g\in G, 
\end{align*}  
where $\Ad(g)_{ij}$ denotes the $ij$-matrix coefficient of  $\Ad(g)$ in the basis $\{X_1, \ldots, \linebreak X_m\}$.

For $\alpha = (\alpha_1, \ldots, \alpha_m)\in \N^m$ with length greater than one the formula is similar, but we have in addition to re-order the derivations, a process that makes polynomial $p_\beta$ in the matrix coefficients appear: 
\begin{align*}
X^\alpha \check{\psi} (g) &= (-1)^{|\alpha|} \big(\,(\Ad(g)X_1)^{\alpha_1} \cdots (\Ad(g)X_m)^{\alpha_m} \psi \,\big) (g^{-1}) \\
&= \sum_{|\beta|\le|\alpha|} \,\big( p_\beta(\Ad(g)_{ij}) X^\beta\psi \,\big)(g^{-1}), \qquad g\in G. 
\end{align*}
 
We want to apply this formula to $\psi=\overline{\phi}\delta$ with $\phi\in \schw_\sigma(G)$. Since the modular function is multiplicative, it is an eigenfunction of the left invariant differential operators $X^\alpha$ (see \cite{varopoulos&al92}). Let  $\lambda_\alpha\in \R$ be such that 
$ X^\alpha\delta=\lambda_\alpha \delta$ ($\alpha\in \N^m$).  
By using Leibniz rule and re-ordering the derivations, we obtain
\[ X^\beta (\bar{\phi} \delta) = \sum_{|\beta_1+\beta_2|\le |\beta|} C_{\beta_1, \beta_2} X^{\beta_1} (\bar{\phi})  X^{\beta_2} \delta = \sum_{|\beta_1+\beta_2|\le |\beta|}  C_{\beta_1, \beta_2} \, \lambda_{\beta_2}\, \delta \, \overline{X^{\beta_1} \phi}. \]
It shows that
\[
X^\alpha (\bar{\phi} \delta)^{\check{ }} (g) =
\delta (g^{-1}) \sum_{|\beta|\le|\alpha|} \,\big( \tilde{p}_\beta(\Ad(g)_{ij}) \overline{X^\beta\phi} \,\big)(g^{-1}), \qquad g\in G, 
\]
with $\tilde{p}_\beta$ polynomials in the matrix coefficients. It follows that 
\begin{align*} 
\|\phi^\ast \|^\infty_{k,\alpha} &= \sup_{g\in G} \left| \sigma^k(g) X^\alpha (\bar{\phi} \delta)^{\check{}}  
(g) \right| \nonumber \\
&\le \sum_{|\beta| \le |\alpha|} \sup_{g\in G} \left| \sigma^k(g) \delta(g^{-1})  \tilde{p}_\beta(\Ad(g)_{ij})  \overline{X^\beta \phi(g^{-1})} \right| \nonumber \\
&= \sum_{|\beta| \le |\alpha|} \sup_{g\in G} \left| \check{\sigma}^k(g) \delta(g)  \tilde{p}_\beta(\Ad(g^{-1})_{ij})  X^\beta \phi(g) \right|\, . 
\end{align*}
Now, on the one hand by Properties~\ref{controle_module} and~\ref{equiv_inverse}, 
\[\check{\sigma}^k(g) \delta(g)\le C \sigma^{rk+m}(g)
, \qquad g\in G.\]
And on the other hand, the operator norm dominates  the modulus of the matrix coefficients, so that there exists $C_\alpha>0$ and  $M_\alpha\in \N$ for which
\[ \left| \tilde{p}_\beta(\Ad(g^{-1})_{ij}) \right| \le C_\alpha \| \Ad(g^{-1}) \|^{M_\alpha}_\opg, \qquad g\in G,\,  |\beta|\le |\alpha|. \]
It follows that 
\begin{equation*} 
\|\phi^\ast \|^\infty_{k,\alpha} \le C_\alpha \sum_{|\beta| \le |\alpha|} \sup_{g\in G} \left| \sigma^{rk+m+M_\alpha}(g)   X^\beta \phi(g) \right|, 
\end{equation*}
which proves the continuity. 
\qed 

\section{More reasons why the function algebra $\SG$ is ``Schwartz-like''
} \label{arbitrary}

\subsection{$L^q$ seminorms} \label{Lqsn}

It is a well known fact from Schwartz theory on Euclidean spaces that $\schw(\R^n)$ can be equipped with families of $L^q$-seminorms, for any choice of $1\le q \le \infty$, and that all these families of seminorms define the same topology on $\schw(\R^n)$. The similar property is true on $\SG$. 

\begin{thm} \label{topologieLq}
Let $1\le q< \infty$. We set  
\[ \| \phi \|^{q}_{k,\alpha}: = \| \sigma^k X^\alpha\phi \|_{L^q(G, dg)}, \qquad \phi\in C^\infty(G)\] with $k\in \N, \ \alpha\in \N^m$.  
The collection of 
$\| \cdot \|^{q}_{k,\alpha}$ 
forms a family of continuous seminorms on $\schw_\sigma(G)$,  which
induces 
the same topology on 
$\schw_\sigma(G)$  than the initial one.  
\end{thm}

\begin{proof} 
We start by showing that the $\| \cdot \|^{q}_{k,\alpha}$ are continuous seminorms on $\SG$. Each $\| \cdot \|^{q}_{k,\alpha}$ being clearly positive homogeneous and subadditive, we only need to prove the continuity. 
This is follows immediately 
from  Property~\ref{compense_volume}, since
\begin{multline*}
\| \phi \|^{q}_{k,\alpha} = \| \sigma^k  X^\alpha\phi \|_{L^q(G,dg)} \, = \, \left( \int_G \frac{1}{\sigma^{qp}(g)} \sigma^{q(k+p)}(g) \left| X^\alpha\phi (g)\right|^{q} dg \right)^{1/q}   \\
\le \, \| \sigma^{k+p}  X^\alpha \phi \|_{L^\infty(G)} \left( \int_G  \frac{1}{\sigma^{qp}(g)} dg \right)^{1/q} \, \le \, C \| \phi \|^{\infty}_{k+p,\alpha}, \qquad \phi \in \SG. 
\end{multline*}
Hence the $\| \cdot \|^{q}_{k,\alpha}$ are continuous seminorms on $\SG$, 

We now show that seminorms $\| \cdot \|^{\infty}_{l,\beta}$ are continuous with respect to the $\| \cdot \|^{q}_{k,\alpha}$, 
which will prove that the $L^\infty$ and $L^q$ seminorms define the same topology on $\SG$. By Hölder inequality and Property~\ref{compense_volume}, it is enough to prove the continuity with respect to the seminorms $\| \cdot \|^{1}_{k,\alpha}$. 

To show the continuity in the $L^1$-seminorms, we need a special case of Sobolev embedding. Since it is a result is of local nature, the Euclidean argument adapts easily. We omit the proof and simply
refer to \cite[Lemma~1.1]{hormander60} for the proof in $\R^m$.
\begin{lem}[Hörmander \cite{hormander60}] \label{sobolev}
There exist a relatively compact open neighbourhood $\Omega$ of identity $e$ in $G$ and a constant $C>0$, such that for every $\phi \in C^\infty(G)$
\[|\phi(e)| \le C \sum_{|\beta|\le m} \int_{\Omega} \left| X^\beta \phi(g) \right| \, dg . \]
\end{lem}
 
Let $\phi \in \SG$. Fix $h\in G$, and denote by $\phi_h$ the translated function $\phi_h(g):= \phi(hg)$. We apply Lemma~\ref{sobolev} to $\phi_h$. Then by left invariance of the $X_j$'s and right invariance of the Haar measure $dg$, we have that 
\begin{align*}
|\phi(h)| 
= |\phi_h(e) | 
&\le C \sum_{|\beta|\le m} \int_\Omega \left|X^\beta \phi_h (g)\right| \, dg  \\
&= C \sum_{|\beta|\le m} \int_\Omega \left|X^\beta \phi(hg)\right| \, dg  \\
&=  C \sum_{|\beta|\le m} \int_{h\Omega} \delta(h^{-1}) \left|X^\beta \phi(g)\right| \, dg, 
\end{align*}
with $C>0$ independent from $h\in G$, $\phi \in \SG$.
 
By properties~\ref{controle_module} and~\ref{equiv_inverse}, $ \delta(h^{-1}) \le C \sigma^{rm}(h)$. 
It follows that \[|\phi(h)| \le C \sigma^{rm}(h)  \sum_{|\beta|\le m} \int_{h\Omega} 
\left| X^\beta \phi(g) \right| dg. \]
Let $k\in \N$ and $\alpha \in \N^m$. Applying the above inequality to $\sigma^k X^\alpha \phi$, and reordering the derivations (using Poincaré-Birkhoff-Witt theorem) yields
\begin{align*} 
|\sigma^k(h) X^\alpha \phi (h)| &\le C \sigma^{k+rm}(h) \sum_{|\beta|\le m}  \int_{h\Omega}  \left| X^\beta X^\alpha \phi(g) \right| \, dg  \\
&\le C \sigma^{k+rm}(h) \sum_{|\gamma|\le m + |\alpha|}   \int_{h\Omega}  \left| X^\gamma\phi(g) \right| \, dg.  
\end{align*}

Let $W_n=\{g\in G; \,  n< \sigma(g) \le n+1\}$. One has $G=\bigcup_{n\in \N} W_n$, with $W_0=\{ e \}$ and $W_n\bigcap W_m = \emptyset$ for $n\neq m$. 
For $n\ge1$, 
\begin{align} \label{supWn}
\sup_{h\in W_n} \left|\sigma^k(h) X^\alpha \phi (h)\right| &\le C \sup_{h\in W_n} \left| \sigma^{k+rm}(h) \sum_{|\gamma|\le m + |\alpha|}   \int_{h\Omega}  \left| X^\gamma\phi(g) \right| \, dg \right| \nonumber \\ 
&\le C (n+1)^{k+rm} \sum_{|\gamma|\le m + |\alpha|}   \int_{W_n\Omega}  \left| X^\gamma\phi(g) \right| \, dg. 
\end{align}
By Property~\ref{sous-polynomial},  for $g=h\omega \in W_n\Omega$ 
\[ \sigma(g)\cdot \sigma(\omega^{-1}) = \sigma(h\omega) \cdot \sigma(\omega^{-1}) \ge C \sigma^{1/s}(h) \ge C n^{1/s}. \]
Since $\Omega$ is relatively compact, $\sigma(\omega^{-1})$ is bounded on $\Omega$, therefore \[\sigma(g)\ge C n^{1/s}. \]
It follows from~(\ref{supWn}) that 
\begin{align*} 
\sup_{h\in W_n} \left|\sigma^k(h) X^\alpha \phi (h)\right| &\le \widetilde{C} n^{k+rm} \sum_{|\gamma|\le m + |\alpha|} \int_{W_n\Omega}  \left| X^\gamma\phi(g) \right| \, dg \\ 
&\le \widetilde{C} \sum_{|\gamma|\le m + |\alpha|} \int_{W_n\Omega} \sigma(g)^{s(k+rm)} \left| X^\gamma\phi(g) \right| \, dg \\
&\le \widetilde{C} \sum_{|\gamma|\le m + |\alpha|} \int_{G} \sigma(g)^{s(k+rm)}\left| X^\gamma\phi(g) \right| \, dg .
\end{align*}
Hence
\begin{align*}
\| \phi \|^\infty_{k, \alpha} 
= \sup_{g\in G} \left|\sigma^k(h) X^\alpha \phi (h)\right|
&\le \widetilde{C} \sum_{|\gamma|\le m + |\alpha|} \int_{G} \sigma(g)^{s(k+rm)}\left| X^\gamma\phi(g) \right| \, dg  \\ 
& =  C \sum_{|\gamma|\le m+|\alpha|} \| \phi \|^1_{s(k+rm), \gamma}, \qquad \phi \in \SG . 
\end{align*}
Therefore the seminorms $\| \cdot \|_{k,\alpha}^\infty$
are continuous with respect to the $\| \cdot \|_{k,\alpha}^1$ ones, which concludes the proof of Theorem~\ref{topologieLq}. 
\end{proof}

The possibility to work with whichever $L^p$-seminorms family is convenient, allows us to prove the following density result.  
 
\begin{thm} \label{DSdense}
The set of smooth functions with compact support $C^\infty_0(G)$ is a dense $\ast$-subalgebra of $\SG$. 
 \end{thm}

\begin{proof}
On the one hand $C^\infty_0(G)\subset \SG$. On the other hand the $\ast$-algebra structure of $\SG$ is inherited from that of $L^1(G,\,dg)$. And finally $C^\infty_0(G)$ is a $\ast$-subalgebra of $L^1(G, \,dg)$. This implies that  $C^\infty_0(G)$ is a $\ast$-subalgebra of $\SG$. 

We are left to show that $C^\infty_0(G)$ is dense in $\SG$. 
We express the density condition in terms of the seminorms $\| \cdot \|^1_{k,\alpha}$ (that define the topology on $\SG$  by Theorem~\ref{topologieLq}): 
If for any given $f\in \SG$, $\epsilon>0$, and finite family $\{(k_i, \alpha_i)\}_{i\in I}\subset \N\times \N^m$, there exists $\phi\in C^\infty_0(G)$ such that $\|f-\phi\|^1_{k_i, \alpha_i}\le \epsilon$ $\forall \, i\in I$, then $C^\infty_0(G)$ is dense in $\SG$. 

Let $f$, $\epsilon$ and $(k_1,\alpha_1), \ldots, (k_q,\alpha_q)$ be as above. 
Consider an approximation to the identity $\{ \rho_j\}_{j\in\N}\subset C^\infty_0(G)$: 
\[ \lim_{j\to\infty} \rho_j \ast \psi = \psi \quad \textrm{in} \quad L^1(G, \,dg), \qquad \psi \in L^1(G,\,dg)  . \]
For simplicity, we assume that the functions $\rho_j$ are all supported inside the compact neighbourhood of identity $U$ defined in \S\ref{notations}. 
We denote by $\mathbf{1}_{U^l}$ the indicator function of the product set $U^l=U\ldots U$. 

We set 
\[f_{j,l} :=  \rho_j \ast(f \cdot \mathbf{1}_{U^l}), \qquad j, l\in \N.\]
Since a function $f_{j,l}$ is the convolution of two compactly supported functions with one of them smooth,  we have $f_{j,l} \in C^\infty_0(G)$. For $j, \, l\in \N$, 
\begin{equation} \label{supp_phij}
\supp f_{j,l} \subset \supp \rho_j \cdot \supp (f \cdot \mathbf{1}_{U^l} ) \subset U \cdot U^l = 
U^{l+1}. \end{equation}
We want to show that there are $j_0,l_0\in \N$ such that
\begin{equation} \label{f-fjl}
\|f-f_{j_0,l_0}\|^1_{k_i, \alpha_i}\le \epsilon,  \qquad  
1\le i \le q, 
\end{equation}
which will conclude the proof of Theorem~\ref{DSdense}.

We write 
\begin{equation} \label{fl+fjl}
\| f - f_{j,l} \|_{k_i,\alpha_i}^1 
\le \| f - f \cdot \indUl \|_{k_i,\alpha_i}^1 + \| f \cdot \indUl - f_{j,l} \|_{k_i,\alpha_i}^1 .
\end{equation}
We will choose $l$ and then $j$ so that each term on the right hand side is smaller than $\epsilon/2$. 

Let us start with the first term. Since the function $\indUl$ is constant almost everywhere on $G$, its derivatives exist and are zero almost everywhere. This implies   
\begin{multline*}
\| f - f \cdot \indUl \|_{k_i,\alpha_i}^1 = \int_G \sigma^{k_i}(g) \left| X^{\alpha_i}\big( f- f\cdot \indUl\big)(g) \right| \, dg  \\
= \int_G \sigma^{k_i}(g) \left| X^{\alpha_i} f(g) - \big(X^{\alpha_i} f\cdot \indUl\big)(g) \right| \, dg
= \int_{G\setminus U^l} \sigma^{k_i}(g) \big| X^{\alpha_i} f(g) \big| \, dg.
\end{multline*}
For all $g\in G\setminus U^l$, we have $l+1 \le |g|_G \le \sigma(g)$. Therefore
\begin{align*}
\int_{G\setminus U^l} \sigma^{k_i}(g) \big| X^{\alpha_i}f(g) \big| \, dg &\le \frac{1}{l+1} \int_{|g|_G \ge l+1} \sigma^{k_i+1}(g) |X^{\alpha_i} f(g)| \, dg  \\
&\le \frac{1}{l+1} \| f \|_{k_i+1,\alpha_i}^1 
\end{align*}

Let $l_0\in \N$ be such that 
\[ \frac{1}{l_0+1} \| f \|_{k_i+1,\alpha_i}^1 \le \epsilon/2, \qquad 1\le i \le q . \]
Then we have 
\begin{equation} \label{ftronquee}
\| f - f \cdot \indUlo \|_{k_i,\alpha_i}^1 <\epsilon/2, \qquad 1\le i \le q.  \end{equation}

We now look at the quantity $\| f \cdot \indUlo - f_{j,l_0} \|_{k_i,\alpha_i}^1$.
By left invariance of the differential operator $X^{\alpha_i}$, we have  
\[X^{\alpha_i} f_{j,l_0} = \rho_j \ast X^{\alpha_i} (f \cdot \indUlo) =  \rho_j \ast ( X^{\alpha_i} f \cdot \indUlo) \]  
almost everywhere on $G$. 
It follows that 
\begin{multline*}
\|  f \cdot \indUlo- f_{j,l_0} \|_{k_i,\alpha_i}^1 = \big\| \sigma^{k_i} X^{\alpha_i}( f \cdot \indUlo - f_{j,l_0}) \big\|_{L^1(G, \,dg)} \\  
= \left\| \sigma^{k_i} \big( X^{\alpha_i} f \cdot \indUlo - \rho_j\ast ( X^{\alpha_i}  f \cdot \indUlo)\big) \right\|_{L^1(G, \,dg)}  \\
= \int_{U^{l_0+1}} \sigma^{k_i}(g) \big| X^{\alpha_i} f(g)  \cdot\indUlo(g) - \rho_j\ast ( X^{\alpha_i}  f \cdot \indUlo)(g)\big| \, dg, 
\end{multline*}
since  $X^{\alpha_i}f \cdot\indUlo$ is supported in $U^{l_0} \subset U^{l_0+1}$, and $\rho_j\ast ( X^{\alpha_i} f \cdot \indUlo)$ is supported in $U^{l_0+1}$ by~(\ref{supp_phij}). The weight function $\sigma$ is bounded on the compact subsets of $G$, so there is $C>0$ such that, for $1\le i \le q$, 
\begin{align*}
\| f \cdot \indUlo- f_{j,l_0} \|_{k_i,\alpha_i}^1  
&\le C \int_{U^{l_0+1}} \big| X^{\alpha_i} f(g)  \cdot\indUlo(g) - \rho_j\ast ( X^{\alpha_i}  f \cdot \indUlo)(g)\big| \, dg \\
&\le C \big\|   X^{\alpha_i} f  \cdot \indUlo - \rho_j\ast ( X^{\alpha_i} f \cdot \indUlo)\big\|_{L^1(G,\,dg)}. 
\end{align*}
Since $\rho_j$ is an approximation identity, we can find  $j_0\in \N$ such  that 
\[\big\|   X^{\alpha_i} f  \cdot \indUlo - \rho_j\ast ( X^{\alpha_i} f \cdot \indUlo)\big\|_{L^1(G,\,dg)} \le \epsilon/2 \]
for all $1\le i\le q$, and so 
\[ \|  f \cdot \indUlo- f_{j_0,l_0} \|_{k_i,\alpha_i}^1 \le \epsilon/2, \qquad 1\le i \le q.\] 

Combining this last estimate with~(\ref{ftronquee}) and~(\ref{fl+fjl}) proves estimate~(\ref{f-fjl}) with $j=j_0$ and $l=l_0$. The Proposition follows. \end{proof}

\subsection{Questioning our conventions}

We discuss the various conventions used to define the function space $\SG$, and explain why $\SG$ is independent of them. 

\subsubsection{The $L^1$ algebra -- does the invariance of the Haar measure play a role in Theorem~\ref{exstnce_pds'}?}

Since there no reason why the Lie group $G$ should be unimodular, it is natural to ask 
whether the result in Theorem~\ref{exstnce_pds'} is typical for the algebra  $L^1(G, dg)$ of functions integrable with respect to the right invariant Haar measure on $G$, or if a similar result holds for $L^1(G, d^lg)$. The answer is that $\SG$ does not distinguish between left and right Haar measure:  

\begin{prop}\label{frechet_alg2}
$\schw_\sigma(G)$ is a dense Fréchet $\ast$-subalgebra of $L^1(G, d^lg)$.  
\end{prop}

\begin{proof}
The proof of Theorem~\ref{exstnce_pds'} adapts almost without changes, the key property being that 
the modular function of $G$ is dominated by a power of $\sigma$ (Property~\ref{controle_module}). Details are left to reader.  
\end{proof}

\subsubsection{The definition of seminorms 
-- what if one considers right invariant derivations?}
 
Proposition~\ref{frechet_alg2} suggests that $\SG$ might be independent  of 
a specific choice of translation invariance. We show below that replacing left invariant derivations by right invariant ones in the seminorms $\| \cdot \|_{k,\alpha}^\infty$, defines the same topology on $\SG$.  

Let $\widetilde{X}_1, \ldots, \widetilde{X}_m$ be the right invariant vector fields on $G$ that agree with the $X_1, \ldots, X_m$ at the origin $e$ of $G$.   
Denote by $\widetilde{X}^\alpha$, $\alpha=(\alpha_1, \ldots, \alpha_m)\in \N^m$, the right invariant differential operators  $\widetilde{X}^\alpha= \widetilde{X}_1^{\alpha_1} \ldots \widetilde{X}_m^{\alpha_m}$ on $G$.

\begin{prop}\label{autres_champs}
Let 
\[ \| \phi \|^{\prime\infty}_{k,\alpha}: = \| \sigma^k \widetilde{X}^\alpha\phi \|_{L^\infty(G)}, \qquad \phi\in C^\infty(G)\] 
with $k\in \N, \ \alpha\in \N^m$.  The collection of $\| \cdot \|^{\prime\infty}_{k,\alpha}$ forms a family of continuous seminorms on $\schw_\sigma(G)$,  which induces the same topology on $\schw_\sigma(G)$  than the initial one.  
\end{prop}

\begin{proof} First, we show that the seminorms $\|\cdot \|^{\prime\infty}_{k,\alpha}$ are continuous on $\SG$. 
The left invariant differentiable operators $X^\alpha$ are related to the $\widetilde{X}^\alpha$'s in the following way: 
\[ \widetilde{X}^\alpha\phi (g) = (-1)^{|\alpha|} X^\alpha \check{\phi}(g^{-1})= (-1)^{|\alpha|} X^\alpha (\overline{\phi^\ast}\delta) (g^{-1}), \qquad \phi \in C^\infty(G). \]
This means that to go from the seminorms relative to right invariant derivations, to the seminorms relative to left invariant derivations (and vice and versa), one needs to handle three operations: involution, multiplication by the modular function, and taking the inverse. By Theorem~\ref{exstnce_pds'} and its proof, these operations are all continuous on $\schw_\sigma(G)$. It implies that each $\| \cdot \|^{\prime\infty}_{k,\alpha}$. 
is  continuous on $\schw_\sigma(G)$. 

To obtain that the seminorms $\| \cdot \|^{\prime\infty}_{k,\alpha}$ define in fact the same topology than the $\| \cdot \|^{\infty}_{k,\alpha}$, we have to show the reverse property of the above one, \emph{i.e.} that the $\| \cdot \|^{\infty}_{k,\alpha}$ are continuous with respect to the $\| \cdot \|^{\prime\infty}_{k,\alpha}$.
The proof relies again on the fact involution, multiplication by the modular function, and taking the inverse are continuous, this time with respect to the $\| \cdot \|^{\prime\infty}_{k,\alpha}$. The computations are similar to those with the seminorms $\| \cdot \|^{\infty}_{k,\alpha}$. We omit the details. \end{proof}

\subsubsection{What if one chooses another subspace $\ce$?} \label{role_de_t}

The question here is whether the particular choice of the subspace $\ce$, used to have an easy expression for the product law of $G$ in real coordinates, has an incidence on the space $\SG$ or not. The following property guarantees that it is not the case. 

\begin{prop} \label{indep_t}
The function space $\SG$ is independent of the choice of the subspace $\ce$.
\end{prop}

\begin{proof} 
Let $\s$ be a subspace of $\g_0$, $\s\neq\ce$,  such that $\s\oplus\n=\g$.  The mapping $I:(s,x) \mapsto (0,x) \cdot s $ is a diffeomorsphism of $\s \times \n$ onto the manifold $G=\g$. For an element $s$ in $ \s$, we denote its coordinates in $\ce\times\n$ by  
$s=(t_s, n_s) = (t_s^1, \ldots t_s^k, n_s)$.    
The coordinate functions $I_j: \s\times\n \to \R$ are given by 
\begin{alignat*}{3} 
&I_j(s,x) =t_s^j,  &
\qquad (s,x)\in s\times n,  &\qquad 1\le j\le k, \\ 
&I_{j+k}(s, x) = P_{j}(n_s,x), &
\qquad (s,x)\in s\times n,  &\qquad 1\le j\le d, 
\end{alignat*}
where the $P_j:\R^d\times\R^d\to \R$ are polynomial maps. 

We endow $\s\times\n$ with a product law similar to~(\ref{pdt_sur_Q}) (turning  $\s\times\n$ into a connected simply connected Lie group isomorphic to $G$). 
Consider a left invariant differential operator of order $r$ on $\s\times\n$. It is mapped by $I$ to a differential operator on $G$, that can be bounded, thanks to computations similar to those from \S\ref{S_algebra},
by a sum of $ \sigma^{l_\alpha}(I(s,x)) X^\alpha$, with $ |\alpha|\le r$ and $l_\alpha\in \N$.

It follows that, if $\tau$ is a weight function on $\s\times\n$  with a definition similar to~(\ref{def_sigma}), the following are equivalent:
\begin{itemize}
\item $f\in C^\infty(\s\times\n)$ and $f$ decrease, with its left derivatives,  
faster than every power of $\tau$.
\item $f\in  C^\infty(G)$ and $f$ decreases, with all its left derivatives,  
faster than every power of $\sigma$.
\end{itemize}
This means that the function space $\SG$ is well defined on $G$, 
independently of the system of real coordinates  considered. 
\end{proof}

\section{Specific cases}
\label{exemples}

We consider three classes of solvable connected simply connected Lie groups. We show how the weight function $\sigma$ captures fundamental geometric properties of the groups in each case, and how the function space 
$\schw_\sigma$ mirrors their specificities.   

\subsection{$G$ nilpotent}  \label{ex_nilp}

The first case that we consider is that of $G$ nilpotent. Then $G$ is equal to its nilradical $N=(\n, \pdtN)$,  so $G=\R^d$ as manifold. The 
the right invariant Haar measure $dg$ is equal to $dn$, the Lebesgue measure on $\R^d$. 

On  connected simply connected nilpotent Lie groups, besides our space of smooth rapidly decreasing functions, there exists a classical notion of Schwartz space that has been introduced and studied long before the present paper (see for instance in \cite{howe77, corwin81}).

\begin{defin}
Let $\schw(G)$ be the space of functions $\phi\in C^\infty(G)$ which are \emph{Schwartz functions as functions on $\R^d$} \textit{i.e.} which satisfy
\begin{equation} \label{top_Rd} 
\sup_{n\in \R^d} \left| |n|^k_{\R^d} D^\alpha \phi(n) \right| < \infty, \qquad k\in \N, \alpha\in \N^d,   \end{equation}
where, for $\alpha=(\alpha_1, \ldots,  \alpha_d)$, $D^\alpha$ denotes the differential operator on $\R^d$ $D^\alpha:= \frac{\partial^{\alpha_1}\ }{\partial n_1^{\alpha_1}} \ldots \frac{\partial^{\alpha_d}\ }{\partial n_d^{\alpha_d}}$. The space $\schw(G)$ is equipped with the topology defined by the  seminorms~(\ref{top_Rd}).
\end{defin}

Our construction allows to recover this classical definition:  

\begin{prop} \label{SG=Snilp} 
The space $\SG$ of smooth functions decreasing $\sigma$-rapidly at infinity on $G$, is equal to  $\schw(G)$ as topological space. 
\end{prop}

\begin{proof}
It is a standard result that, because the product on $G=(\n, \pdtN)$ is given by a polynomial mapping in the coordinates on $\R^d$, the left invariant vector fields $X_j$ ($1\le j\le d$) have  expressions of the form:
\begin{equation} \label{X=sumD}
\left({X_j}\right)_n  = \sum_{i=1}^d P^j_i(n) \left(\frac{\partial }{\partial n_i}\right)_n, \qquad n\in G, \end{equation}
where the $P^j_i:\R^d\to \R$ are polynomial maps, and reciprocally, that the coordinate derivatives $\frac{\partial \ }{\partial n_j}$ can be written as
\begin{equation} \label{D=sumX} 
\left( \frac{\partial}{\partial n_j}\right)_n  = \sum_{i=1}^d Q^j_i(n) \left(X_i\right)_n, \qquad n\in G,  \end{equation}
where $Q_j^i$ is a polynomial map from $\R^d$ to $\R$ for every $ i, j\in \{1, \ldots, d\}$ (see for instance \cite{ludwig-molitorbraun95}). 

The weight function $\sigma$ is comparable to powers of the Euclidean norm on $\R^d$. Indeed, by definition, 
\[ \sigma(n) = \max\big( \|\Ad (n)\|_\opn, \,  \|\Ad(-n)\|_\opn \big) \cdot \big(1+2\,|n|_N\big) 
, \qquad n\in G. \] 
So, by~(\ref{controle_u}) and~(\ref{domination_Ad}), there are $C>1$ and $q, q'\in \N^\ast$, such that
\begin{equation} \label{poids_nilp}
C^{-1} |n|_{\R^d}^{1/q} \le \sigma(n) 
\le C(|n|_{\R^d}+1)^{q'+1}, \qquad n\in G.  
\end{equation}

Let us show that 
\begin{equation*} 
\schw(G) \subset \SG,  
\end{equation*}
with continuous inclusion map. 
Let $\phi$ be in $\schw(G)$. 
For every $k\in \N$ and $\alpha \in \N^d$, we have by~(\ref{poids_nilp}) that
\begin{equation} \label{sens1}
\sigma^k(n) \left| X^\alpha \phi(n) \right|  \le C^k(|n|_{\R^d}+1)^{k(q'+1)}
 \left| X^\alpha \phi(n) \right|, \quad n\in G .
\end{equation} 
By~(\ref{X=sumD}), there are $C_\alpha>0$ and $k_\alpha\in \N$ such that 
\[\left| X^\alpha \phi (n) \right| \le C_\alpha (1+|n|_{\R^d})^{k_\alpha} \sum_{|\beta| \le |\alpha|} \left| D^\beta \phi(n)\right|,  \qquad n \in G, \ \phi \in \schw(G) .\] 
We inject the estimate above in~(\ref{sens1}) and take the supremum. It  shows
\begin{equation*} 
\sup_{n\in G} \left| \sigma^k(n)  X^\alpha \phi(n) \right| 
\le C_\alpha' \sup_{n\in G} \Big( (1+|n|_{\R^d})^{k'_\alpha} \sum_{|\beta|\le |\alpha|} \left| D^\beta \phi(n)\right|\Big) .  
\end{equation*} 
with $C_\alpha'$ independent of $\phi$. This shows that $\schw(G)$ is contained into $\SG$, with continuous inclusion map.  

A similar argument proves that $\SG \subset \schw(G)$ with continuous inclusion map, so the two spaces coincide and have the same topology. 
\end{proof}

\subsection{$G$ exponential solvable} \label{ex_exponentiel}

In this second case, we assume that 
the solvable Lie group $G$ is  \emph{exponential}, which means that the exponential map $\exp:\g\to G$ is a diffeomorphism. To avoid falling into the scope of the previous example, we suppose furthermore that $G$ is not nilpotent.  
In this situation, $\ce$ has dimension $k\ge 1$ and 
$G=\ce\times\n=\R^k\times\R^d$ as manifold.   
Note that, since $G$ is exponential, the operators $\ad t$ with $t\in \ce$, have no purely imaginary eigenvalues (\cite[Theorem~I.2.1]{bernat&al72}). 

On such groups, no natural definition of Schwartz space arises as clearly as in the situation described in \S\ref{ex_nilp}. We compare $\SG$ with  $\mathcal{ES}(G)$, a space of smooth rapidly decreasing functions introduced by J. Ludwig in \cite{ludwig83}. In the  particular case where $\ce$ is one-dimensional, we show  that the two function spaces agree with one another.

\begin{defin}
Let $\mathcal{ES}(G)$ be the space of functions $\phi\in C^\infty(G)$ such that
\begin{equation} \label{top_ludwig}
 \sup_{(t,n)\in \R^k\times\R^d} \left| e^{r | t |_{\R^k}} | n |_{\R^d}^j D^{\alpha} \phi(t,n) \right| < \infty, \qquad r\ge 0, \, j\in \N, \, \alpha \in \N^{k+d},
\end{equation}
where, for $\alpha=(\alpha_1, \ldots, \alpha_{k+d}) \in \N^{k+d}$, $D^\alpha$ denotes the differential  operator  on  $\R^k\times\R^d$  $D^{\alpha}:= \frac{\partial^{\alpha_1}\ }{\partial t_1^{\alpha_1}} \ldots \frac{\partial^{\alpha_{k+d}\ }}{\partial n_d^{\alpha_{k+d}}}$. The space $\schw(G)$ is equipped with the topology defined by the  seminorms~(\ref{top_ludwig}).
\end{defin}

\begin{prop} \label{SG>SE}
The space $\mathcal{ES}(G)$ is densely contained in the Fréchet algebra $\SG$ of smooth functions decreasing $\sigma$-rapidly at infinity on $G$, with continuous inclusion map. 

In the case where $\ce$ has dimension one, $\mathcal{ES}(G)$ and $\SG$ are equal as topological spaces. 
\end{prop}

\begin{rem}
The exponential solvable Lie groups $G$ for which $\dim \ce=1$ are rank one $NA$ groups. A typical example is the $ax+b$ group, which is a semidirect product of $\R$ with $\R$. 
\end{rem}

\begin{proof} 
We start by showing that $\mathcal{ES}(G)\subset \SG$. Our first task is prove that 
there are $C>0$ and $l\in \N$,  such that 
\begin{equation} \label{sigma<}
\sigma(g) \le C e^{C | t |_{\R^k}} (1+| n |_{\R^d})^{l}, \qquad g=(t,n)\in G.\end{equation} 
By Property~\ref{sous-polynomial}, $\sigma$ is sub-polynomial. Writing $(t,n)=(0,n)\cdot(t,0)$ shows that it is enough to prove (\ref{sigma<}) for $g=(t,0)$ and for $g=(0,n)$.

Let $t\in \ce$. By definition, 
\[ \sigma(t,0)= \max\big( \|\Ad (t,0)\|_\opg, \,  \|\Ad(-t,0)\|_\opg \big) \cdot \big(1+|(t,0)|_G\big). \]
We estimate the two factors independently. For $|(t,0)|_G$, we have 
\begin{equation} \label{tG/tRk}
|(t,0)|_G \le C \big( | t |_{\R^k} + 1 \big), \qquad t\in \ce.  
\end{equation}
Indeed, assume that $U \supset \left\{(t, 0)\in \ce\times\n; \  | t |_{\R^k} \le 1\right\}. $
Then $(t,0)$ can be written as product of $ [ | t |_{\R^k} ] +1$ elements of $U$:
\[ (t, 0) = \left(\frac{t}{ [| t |_{\R^k}]  +1} \, , \, 0 \right)  \, \cdots \,  \left(\frac{t}{[ | t |_{\R^k} ]  +1} \, , \, 0 \right),   
\]
which implies~(\ref{tG/tRk}) with $C=1$. The general case where $U$ does not necessarily contain 
$\left\{(t, 0)\in \ce\times\n; \  | t |_{\R^k} \le 1\right\}$ follows from~(\ref{equiv-lgueurs}). 

Now we turn to the first factor in the expression of $\sigma(t,0)$. The operator $\ad$ is linear continuous on $\g$, so there is $C>0$ such that, for every $ (t,0) \in \ce\times\n$, 
\[
\| \Ad(t,0) \|_\opg= \|e^{\ad(t,0)} \|_\opg \le e^{\|\ad(t,0)\|_\opg}  \le e^{C |(t,0)|_{\R^m}} = e^{C |t|_{\R^k}}.
\]
Combined with~(\ref{tG/tRk}), this shows that 
\[ \sigma(t,0) \le C e^{C |t|_{\R^k}} , \qquad t\in\ce \]
which proves~(\ref{sigma<}) for $g=(t,0)$. 

Let us now prove~(\ref{sigma<}) for $n\in \n$.
We have \[\sigma(0,n) = \max\big( \|\Ad (0,n)\|_\opg, \,  \|\Ad(0, -n )\|_\opg \big) \cdot \big(1+|(0,n)|_G + |n|_N \big) . \]
The right hand side factor is easy to bound. By~(\ref{comparaison_QN}) and~(\ref{controle_u}), 
\begin{equation} \label{poursigma_n}
|(0,n)|_G + |n|_N \le 2 |n|_N  \le C ( 1+ |n|_{\R^d} ), \qquad n \in \n.
\end{equation}
To estimate the left hand side factor, recall that $\n$ is the nilradical of $\g$, so that $\ad{(0,n)}^{m}=0$ for every $n\in \n$.  Therefore 
\begin{multline*} 
\| \Ad(0,n) \|_\opg = \| e^{\ad(0,n)} \|_\opg =  \left\| \sum_{j=0}^m \frac{ \ad(0,n)^j }{j!} \right\|_{\opg} \\ 
\le  \sum_{j=0}^m \frac{\| \ad(0,n) \|_{\opg}^j}{j!} \le C \sum_{j=0}^m \frac{| (0,n) |_{\R^m}^j }{j!} \le C (1+|n|_{\R^d})^m, \qquad n \in \n. 
\end{multline*}
Together with~(\ref{poursigma_n}), this shows that \[\sigma(0,n) \le C (1+|n|_{\R^d})^{m+1}, \qquad n \in \n.  \] This completes the proof of estimate~(\ref{sigma<}).

Our second task is to compare the differential operators $X^\alpha$ with the $D^\alpha$. 
For this, it will be convenient to consider a distinguished basis of the Lie algebra $\g$. 
Let $\{e_1, \ldots e_d\}$ be the canonical basis of $\R^d$ and $\{e_{d+1}, \ldots e_{d+k}\}$ the canonical basis of $\R^k$. As elements of $\g=\ce\times\n$, the $e_j$'s induce left invariant vector fields on $G$ by
\[X_j\phi(g)= \left. \frac{\partial}{\partial t} \phi(g \exp_G(t e_j))\right|_{t=0},  \qquad \phi \in C^\infty(G), \]
that are linearly independent. We take $\{X_1, \ldots, X_{d+k}\}$ as basis of the Lie algebra $\g$.   

For $1\le j\le d$, the $e_j$'s induce also left invariant vector fields on $N$ 
\[ \widetilde{X}_j\psi(n) = \frac{\partial}{\partial t} \psi (n \pdtN  \exp_N(t e_{j})), \qquad \phi\in C^\infty(N). \] 

The vector fields $X_j$ are related to the $\widetilde{X}_i$ in the following way:
\begin{align}
\left( X_j \right)_{(t,n)} &= e^{\ad t} \left( \widetilde{X}_j\right)_n, \qquad 
 1\le j\le d \label{X-Xtilde1} 
\\
\left( X_{d+j} \right)_{(t,n) } &= \left(\frac{\partial}{\partial t_j}\right)_{t} + \sum_{i=1}^d R_i^j(t) \left( \widetilde{X}_i \right)_n, \qquad 
 1\le j \le k, \label{X-Xtilde2}
\end{align}
where the $R_i^j:\R^k\to \R$ are polynomial map. 

On the one hand, the terms $ \|e^{\ad t }\|_{\opn}$ and  $ |R_i^j(t)|$ ($1\le i,j\le d$) are dominated  by $C e^{C |t|_{\R^k}}$ for $C>0$ large enough. On the other hand, we know from the previous example (\S\ref{ex_nilp}), that for $\alpha\in \N^d$ there are $C_\alpha>0$ and $k_\alpha \in \N$ such that 
\[\left| \widetilde{X}^\alpha \psi (n) \right| \le C_\alpha (1+|n|_{\R^d})^{k_\alpha} \sum_{|\beta| \le |\alpha|} \left| D^\beta \psi(n)\right|,  \qquad n \in N .\] 
It follows that for every $\alpha\in \N^{d+k}$, there are $C_\alpha>0$ and $k_\alpha\in \N$ such that, for all $(t,n) \in G$,
\begin{equation} \label{Xtn<Dtn}
\left| X^\alpha \phi (t,n) \right| \le C_\alpha (1+|n|_{\R^d})^{k_\alpha} e^{C_\alpha |t|_{\R^k}} \sum_{|\beta| \le |\alpha|} \left| D^\beta \psi(t,n)\right|.   
\end{equation}
 
Combining the estimates~(\ref{sigma<}) and~(\ref{Xtn<Dtn}) and taking the supremum, we obtain that for every $l\in \N$ and $\alpha\in \N^{d+k}$, there are $C>0$ and $j\in \N$ such that 
\[ \sup_{g\in G} \left| \sigma^k(g) X^\alpha \phi(g) \right| \le C \sup_{(t,n)\in G} 
 (1+|n|_{\R^d})^j e^{C |t|_{\R^k}} 
\sum_{|\beta| \le |\alpha|} \left| D^\beta \psi(t,n)\right| 
.\] 
It proves that $\mathcal{ES}(G)$ is contained in $\SG$,  and that the inclusion map is continuous. The density comes from the obvious fact that $C^\infty_0(G) \subset \mathcal{ES}(G)$. 

In the special case where $\ce=\R$, the key observation is that, for some $\mu>0$,
\[ e^{\mu |t|_{\R}} \le \max( \|e^{\ad t} \|_\opn, \|e^{- \ad t}\|_\opn ) , \qquad t\in \ce  \]
(which is not necessarily true when $k>1$).
By using arguments similar to those in \S\ref{pf_1'}, one obtains easily that  
\[ e^{|t|_\R} |n|_{\R^d} \le C \sigma^j(t,n), \qquad (t,n)\in G . \]
The domination of the differential operators $D^\alpha$ ($\alpha\in \N^{d+1}$) by sums of $\sigma^{j}
X^\beta$ follows from the formulas~(\ref{X-Xtilde1}),~(\ref{X-Xtilde2}) and estimate~(\ref{D=sumX}).  

With these controls on the weight and on the differential operators, one can dominate  each seminorm on $\mathcal{ES}(G)$ by a finite sum of $\| \cdot \|_{k,\alpha}^\infty$. This proves the reverse inclusion: $\SG\subset \mathcal{ES}(G)$, with the equivalence between the topologies. \end{proof}

\subsection{From the motion group of $\R^2$  
to groups with polynomial growth}

Take $G=M(2)$, the universal covering group of the group of Euclidean motions of the plane. At the level of the Lie algebra $\g$, there is basis $\{X_1,X_2,X_3\}$  which satisfies the commutation relations  
\[ [X_1,X_2]=X_3, \qquad [X_1,X_3]=-X_2, \qquad [X_2,X_3]=0   . \]
The nilradical $\n=\vect(X_2,X_3)$ is commutative, and we take $\ce = \vect(X_1)$ (see the discussion in \S\ref{role_de_t}).
At the group level, the product law~(\ref{pdt_sur_Q}), in coordinates $(t,n_1,n_2)\in \R^3$, is 
\[ (t,n_1,n_2)\cdot (t',n_1',n_2') = (t+t',\, n_1+   
n_1'\cos t + n_2'\sin t, n_2 - n_1'\sin t + n_2'\cos t) .\]

Note that this example is not covered by  \S\ref{ex_nilp} and \S\ref{ex_exponentiel}, because $G$ is a (non-nilpotent) solvable Lie group which is not exponential (because $\ce$ acts on $\n$ by rotations).  

The weight $\sigma$ is easy to estimate here, all the terms involved in its expression being more or less explicit. We have, for every $g=(t,n_1,n_2)\in G$,
\begin{equation} \label{sigma-M2}
C^{-1} \left(1+|(t,n_1,n_2)|_{\R^3}\right)\le \sigma(g) \le C  \left(1+|(t,n_1,n_2)|_{\R^3}\right)^2, 
\end{equation}
with $C>1$ independent of $g$. 

We denote again by $X_1,X_2,X_3$ the left invariant vector fields induced by the basis of the Lie algebra. Their expression in coordinates $(t,n_1,n_2)$ is 
\[ X_1 = \frac{\partial}{\partial t}, 
\qquad X_2 = \cos t  \frac{\partial}{\partial n_1} - \sin t \frac{\partial}{\partial n_2}, \qquad X_3 = \sin t  \frac{\partial}{\partial n_1} + \cos t 
\frac{\partial}{\partial n_2} . \]
This implies in return that 
\[ \frac{\partial}{\partial n_1} = \cos t X_2 + \sin t X_3, \qquad \frac{\partial}{\partial n_2} = - \sin t X_2 + \cos t X_3 .\] 
These relations allow to dominate 
every differential operator with constant coefficients in $\R^3$ by sums of $X^\alpha$ ($\alpha\in \N^3$) and vice and versa.  
By the same reasoning than that used to prove Propositions~\ref{SG=Snilp} and~\ref{SG>SE}, this yields, together with the estimate of $\sigma$ (\ref{sigma-M2}), the following result:   

\begin{prop} \label{M2}
The space $\SG$ of smooth functions decreasing $\sigma$-rapidly at infinity on $G=M(2)$ is equal, as topological space, to the usual Schwartz space on $\R^3$.
\end{prop}

With respect to the function space $\schw_\sigma$, the motion group $M(2)$ can be seen as a model case for non-nilpotent solvable Lie groups with polynomial growth of the volume. Proposition~\ref{M2} is a special case of the more general result: 

\begin{prop}
Let $G$ be an  $m$-dimensional connected simply connected solvable Lie group with polynomial growth of the volume. Then the Fréchet space $\SG$ of smooth functions decreasing $\sigma$-rapidly at infinity on $G$ is the usual Schwartz space on $\R^m$.  
\end{prop}

\begin{proof}
If $G$ is nilpotent, the statement is true by Proposition~\ref{SG=Snilp}. If $G$ is non-nilpotent, the subspace $\ce$ is non-zero and the operators $\ad t$ ($t\in\ce$) have all their eigenvalues purely imaginary. 
The ingredients in the proof of 
Proposition~\ref{M2}, \textit{i.e.}, on the one hand the upper and lower estimates of $\sigma$ by integer powers of the Euclidean norm, and on the other hand the control of the differential operators with constant coefficients on $\R^m$ by left invariant differential operators on $G$ and reciprocally,  extend in an easy way to the general setting, leading to the result. We omit the details.  
\end{proof}

\section{Tempered distributions} \label{T_distributions}

We leave examples, and come back to the general case. The notations are those from \S\ref{preliminaires} and \S\ref{S_algebra}. 

\begin{defin}
The topological dual space $\schw'_\sigma(G)$ of $\schw_\sigma(G)$, is called the space of $\sigma$-tempered distributions on $G$. 
For $T\in \schw'_\sigma(G)$ and $\phi\in \schw_\sigma(G)$, we denote the evaluation of $T$ on $\phi$ by $<T,\phi>$. 
\end{defin}

The terminology \emph{tempered} is justified by the fact that in our setting, like in the Euclidean setting,  there is a certain growth restriction at infinity on the elements of $\schw'_\sigma(G)$, as  will be seen now.   

\begin{defin}
We say that a measurable function $f$ on $G$ is $\sigma$-slowly increasing at infinity, if there is $k\in \N$ such that \[\| \sigma^{-k} f \|_{L^\infty(G)} <\infty. \] 
\end{defin}

It is a standard verification that measurable functions $f$ increasing $\sigma$-slowly at infinity on $G$, determine tempered distributions $\big[f\big]$, via the correspondence:
\[ \big[f\big]: \, \phi \mapsto \int_G \phi(g) f(g) dg . \]
In particular, this is the case of continuous functions increasing $\sigma$-slowly at infinity 
on $G$. As a matter of fact, these functions, together with their distributional derivatives, exhaust the space $\schw'_\sigma(G)$. 

Before giving a precise statement of this result, let us define the notion of distributional derivative. 

\begin{defin}
Let $T\in \schw'_\sigma(G)$ and $\alpha \in \N^m$. The mapping  
\[
X^\alpha T : \, \phi \mapsto (-1)^{|\alpha|} <T, X^\alpha \phi> 
\]
defines a continuous linear functional on $\SG$, \textit{i.e.} a $\sigma$-tempered distribution on $G$. We say that $X^\alpha T$ is a \emph{distributional derivative} of $T$.\end{defin}
 
Note that if a function $f$ is differentiable up to the order $|\alpha|$ in the classical sense, and if it has all its derivative $\sigma$-slowly increasing, we have $X^\alpha \big[ f \big] = \big[ X^\alpha f \big] $ (follows easily from the right invariance of the Haar measure $dg$).

\begin{thm} \label{D-temperee} Let $T\in \schw'_\sigma(G)$. There exists  $M\in \N$ and 
a finite family $\{f_\alpha; \, \alpha \in \N^m, \,  |\alpha|\le M \}$  
of continuous functions increasing $\sigma$-slowly at infinity on $G$,  
such that 
\[ T = \sum_{|\alpha| \le  M } X^\alpha \big[ f_\alpha \big] .\]
\end{thm}

\begin{proof} We follow the lines of the proof of the classical result in $\R^n$ \cite[Theorem 51.6]{treves67}. 
Let $T\in \schw'_\sigma(G)$. By Theorem~\ref{topologieLq}, $T$ is continuous with respect to family of seminorms $ \| \phi \|_{j,\alpha}^1 = \| \sigma^j X^\alpha\phi \|_{L^1(G,\, dg)}$ on $\SG$. So 
there exist $C>0$ and $j\in \N$ such that 
\begin{equation}  \label{cont_T}
 \left| <T,\phi> \right| \le 
C \sum_{
|\alpha| \le  j } 
\| \sigma^j X^\alpha\phi \|_{L^1(G, dg)},
\qquad \phi \in \SG. \end{equation}
Let $A$ be the number of multi-indices $\alpha \in \N^m$ such that $|\alpha|\le j$. 
Define \begin{align*} 
J: \quad \SG  
&\to \underbrace{L^1(G,\, dg) \times \ldots \times L^1(G,\, dg)}_{A \textrm{ times}} =  L^1(G,\, dg)^A  \\ 
\phi &\mapsto \left(\sigma^j X^\alpha \phi \right)_{|\alpha|\le j}   .
       \end{align*}
The map $J$ is linear and, because $\sigma $ is strictly positive on $G$, it is injective. 

Formulated in terms of $J$,
estimate~(\ref{cont_T}) means that the linear functional $J \phi \mapsto <T,\phi> $ is continuous on $J \SG$ 
equipped with the  topology inherited from $L^1(G,\, dg)^A$. By Hahn-Banach Theorem, it can be extended to a continuous linear functional on the whole space $L^1(G,\, dg)^A$. Since the dual space of $L^1(G,\, dg)^A$ is $L^\infty(G,\, dg)^A$, we have in particular that there are $A$ functions $h_\alpha$ ($|\alpha|\le j$) in $L^\infty(G,\, dg)$ 
such that 
\[ <T,\phi> = \sum_{|\alpha|\le j } \int_G h_\alpha(g) \sigma^{j}(g) X^\alpha \phi(g) dg, \qquad \phi\in \SG,  
\]
which means 
\begin{equation} \label{T_Linfini} 
T= \sum_{|\alpha|\le j } (-1)^{|\alpha|} X^\alpha \big[ \sigma^j h_\alpha \big] .  \end{equation}

At this stage, we have written $T$ as a sum of distributional derivatives of measurable functions $\sigma$-slowly increasing at infinity on $G$. The next step is to have a sum of derivatives of continuous functions on $G$.

Let $g_0=(t_0,n_0) \in G$ with $t_0=(t_1, \ldots, t_k)\in \ce 
$ and $n_0=(n_1, \ldots, n_d) \in \n $. 
For each $|\alpha|\le j$, we set  
\[ {h_\alpha^\flat}(g)= \int_0^{t_1}\!\!\!\ldots \int_0^{t_k} \!\!\!\int_0^{n_1}\!\!\!\ldots \int_0^{n_d} \big(\sigma^j h_\alpha\big)(t,n) \,dt \, dn .\] 

Recall that the right invariant Haar measure $dg=dt\,dn$ is the Lebesgue measure on $\R^m=\R^k\times\R^d$. 
Since $\sigma$ is continuous almost everywhere on $G$ and $h_\alpha$ is in $L^\infty(G, \, dg)$, then the function $\sigma^j h_\alpha$ is locally integrable on $G$ with respect to $dg$, equivalently on $\R^m$ with respect to the Lebesgue measure. It follows that  
$h_\alpha^\flat$ is continuous as function on $\R^m$ 
(hence on $G$). Furthermore, if we denote by $D$ the differential operator $ D:= \frac{\partial\ }{\partial t_1} \ldots \frac{\partial\ }{\partial t_k}\frac{\partial\ }{\partial n_1} \ldots \frac{\partial\ }{\partial n_d}$ on $\R^m$, we have that, for all $\phi\in C_0^\infty(\R^m)$, 
\begin{equation} \label{Dh.}
\int_{\R^m} h_\alpha^\flat(g) D\phi(g) dg = (-1)^m \int_{\R^m} \sigma^j(g) h_\alpha(g)  \phi(g) dg. \end{equation}

By general theory, the first order derivative in $\R^m$ are related to the left invariant vector fields on $G$ by the mean of smooth functions $\tau^i_l$, $\nu^i_l$, in the following way: 
\[ \frac{\partial\ }{\partial t_i} = \sum_{l=1}^m \tau^i_l X_l \quad (1\le i\le j),  \qquad\textrm{and} \qquad \frac{\partial\ }{\partial n_i} = \sum_{l=1}^m \nu^i_l X_l \quad (1\le i\le d). \]
It follows that there are functions $\eta_\beta \in C^\infty(G)$ ($\beta\in \N^m, |\beta| \le m$), such that 
\begin{equation*} \label{DX} 
D = \sum_{|\beta|\le m} \eta_\beta X^\beta  
\end{equation*}

Because the notions of smooth functions on $G$ and on $\R^m$ coincide (and similarly for $C^\infty_0(G)$ and  $C^\infty_0(\R^m)$), we can inject this equality  in~(\ref{Dh.}). We obtain that, for all $\phi \in C^\infty_0(G)$, 
\[ \sum_{|\beta| \le m}\int_{G} h_\alpha^\flat(g) \eta_\beta X^\beta \phi(g) dg = (-1)^m \int_{G} \sigma^j(g) h_\alpha(g)  \phi(g) dg. \]
By density of $C^\infty_0(G)$ in $\SG$, the equality is valid for all $\phi\in\SG$, which means precisely that 
\[ \sum_{|\beta| \le m} X^\beta  \left[(-1)^{|\beta|+m} \eta_\beta h_\alpha^\flat \right] = \big[ \sigma^j h_\alpha \big] ,\]
and the functions $\eta_\beta h_\alpha^\flat$ ($|\beta|\le m$) are continuous on $G$. 

We apply this identity to~(\ref{T_Linfini}), then write every differential monomial $X^\alpha X^\beta$ as a finite sum
\[ X^\alpha X^\beta = \sum_{i\in I} X^{\gamma_i}, \qquad |\gamma_i|\le|\alpha| + |\beta|,  \] 
and we gather the terms with same order of derivation. We obtain a collection of continuous functions $f_\alpha$ which verify 
\[ T = \sum_{|\alpha|\le j} \sum_{|\beta|\le m} X^\alpha X^\beta \left[ (-1)^{|\beta|+|\alpha|+m} \eta_\beta h_\alpha^\flat\right] = \sum_{|\gamma|\le j+2m} X^\gamma \big[f_\gamma \big]. \] This gives the Theorem with $M=j+2m$. 
\end{proof}

\section{ Rapidly decreasing functions on product  group and Kernel Theorem} 
\label{S_thm}

Let $G_1$, $G_2$ be two connected simply connected solvable Lie groups with  respective dimensions $m_1$ and $m_2$. We denote by $\g_1$ and $\g_2$ their respective Lie algebras. 
For $i=1,2$, let $\n_i$ the nilradical of $\g_i$.  
By~\S\ref{decomp_G}, there are subspaces $\ce_i\subset \g_i$, and polynomial mappings $P_i:\ce_i\times\ce_i\to \n_i$ such that $G_i=\ce_i\times\n_i$ with product law in coordinates $(t_i,n_i)\in \ce_i\times\n_i$
\[(t_i,n_i)\cdot(t_i',n_i')=(t_i+t_i',P(t_i,t_i') \pdtN n_i \pdtN  e^{\ad t_i}n_i'). \]

In this section, we are interested in Schwartz functions on the product group $G=G_1\times G_2$.  
In order to define these functions,  we need to write the Lie group $G$ 
in a way similar to that in \S\ref{decomp_G}.   

Let $\g$ denote the Lie algebra of $G$. It decomposes as $\g=\g_1\oplus\g_2$. The nilradical of $\g$ is then $\n=\n_1\oplus\n_2$. To write the product law in $G$ in real coordinates, we must choose, inside a certain nilpotent subalgebra of $\g$, 
a complementary subspace $\ce$ to the nilradical $\n$ in $\g$. 
The important point is that, inside this given nilpotent subalgebra, we are free to choose $\ce$ as we like (see Proposition~\ref{indep_t}). The idea is to choose $\ce$ coherent with respect to the product group structure. 
  
Assume that the subspaces $\ce_i\subset \g_i$ are defined with respect to elements $X_i\in \g_i$ in general position relatively to the roots of $\g_i$ ($i=1,2$). 
It is easy to check that $X=(X_1,X_2)\in \g$ is in general position with respect to the set of roots of $\g$. With the notations from \S\ref{decomp_G}, we have then 
\begin{align*}
\g_{\C,0}
&= \big\{ Y\in \g_\C;\, \ad(X)^{\dim \g}(Y)=0 \big \} \\    &= \Big\{ (Y_1,Y_2)\in {\g_1}_\C \oplus {\g_2}_\C;\, \big(\ad(X_1)^{\dim \g}(Y_1),\ad(X_2)^{\dim \g}(Y_2)\big) =0 \Big\}  \\ 
&= {\g_1}_{\C, 0} \oplus   {\g_2}_{\C, 0}.                 \end{align*}
Hence 
\[ \g_0 = \g_{\C, 0}\cap \g = \big({\g_1}_{\C, 0}\cap \g_1 \big) \oplus  \big({\g_2}_{\C, 0} \cap \g_2) =  {\g_1}_0 \oplus   {\g_2}_0 . \]

The complementary subspace $\ce$ 
to $\n$ in $\g$, has to be chosen inside the subspace $\g_0$. We take $\ce=\ce_1\oplus\ce_2 \subset  {\g_1}_0 \oplus   {\g_2}_0$. Then we have 
$G=\ce\times\n$ with product law 
\[(t,n)\cdot(t',n')=(t+t',P(t,t') \, \cdot_\textrm{\tiny CBH} n \,
\cdot_\textrm{\tiny CBH} e^{\ad t}n'), 
\qquad (t,n), (t',n')\in \ce\times\n, \]
where  $P:\ce\times\ce\to \n$ 
is the polynomial mapping defined by 
\[  \begin{array}[t]{cccc} 
P: &(\ce_1\oplus\ce_2) \times (\ce_1\oplus\ce_2) 
& \longrightarrow &  
\n_1\oplus \n_2 \\ 
& \big((t_1,t_2),\,(t'_1,t'_2)\big) &\longmapsto & \big(P_1(t_1, t'_1), \, P_2(t_2, t'_2) \big).
   \end{array}
\] 

We are now ready to define the weight functions that will characterize the spaces of rapidly decreasing functions on $G_1$, $G_2$ and $G$. 

For $i=1,2$, let $U_i$ be a symmetric compact neighbourhood of identity in $G_i$. The product set $U=U_1\times U_2$ is a symmetric compact neighbourhood of identity in $G$. We denote by $|\cdot|_{G_i}$, resp. $|\cdot|_G$, the length of an element of $G_i$, resp. $G$,  with respect to $U_i$, resp. $U$. 
We have \[ |g|_G=\max\big(|g_1|_{G_1},\, |g_2|_{G_2}\big), \qquad g=(g_1,g_2)\in G.  \] 

Similarly, if $|\cdot|_{N_i}$ ($|\cdot|_N$) is the length of an element of the nilradical $N_i$ of $G_i$ (of the nilradical $N$ of $G$)  with respect to $U_i\cap N_i$ ($U\cap N$), we have \[ |n|_N=\max\big(|n_1|_{N_1},\, |n_2|_{N_2}\big), \qquad n=(n_1,n_2)\in N.  \] 

We set 
\begin{align*}
\sigma_i(g) &= \max\big( \|\Ad (g)\|_\opgi, \,  \|\Ad(g^{-1})\|_\opgi \big) \cdot \big(1+|g|_{G_i}+|n|_{N_i}\big) \quad (i=1,2)
\\  \sigma(g) &= \max\big( \|\Ad (g)\|_\opg, \,  \|\Ad(g^{-1})\|_\opg \big) \cdot \big(1+|g|_G+|n|_N\big),\end{align*}
where $\|\cdot\|_\opgi$ ($\|\cdot\|_\opg$) denotes the operator norm on $\g_i$ ($\g$) equipped with the Euclidean norm.

\begin{thm} \label{Stensor}
Let the Lie groups $G_1$, $G_2$, $G$, with the respective weight functions $\sigma_1$, $\sigma_2$, $\sigma$, be as above.
Then $\schw_\sigma(G)$ is isomorphic as Fréchet space to the projective completion $ \schw_{\sigma_1}(G_1) \widehat{\otimes}_\pi \schw_{\sigma_2}(G_2).$
\end{thm}

In order to prove the result, we need to show first that the Fréchet spaces of rapidly decreasing functions defined in this paper, are nuclear spaces (we refer to  \cite{treves67} for the general theory of nuclear topological vector spaces).  

\begin{prop} \label{Snuclear}
Let $G$ and $\sigma$ be as in \S\ref{preliminaires} and \S\ref{S_algebra}. The  
space $\SG$ of smooth functions decreasing $\sigma$-rapidly at infinity on $G$ is nuclear.  
\end{prop}

\begin{proof}[Proof of Proposition~\ref{Snuclear}] 
We use a criterion for nuclearity 
due to L. Schweitzer \cite{schweitzer93a} (the statement given here is slightly weaker than the one in~\cite{schweitzer93a}, to avoid the use of 
new definitions). 

\begin{lem}[Schweitzer {\cite[Thm 6.24]{schweitzer93a}}] \label{schweitzer}
Let $H$ be a Lie group and $\tau$ be weight function on $H$. Assume that $\tau$ is sub-polynomial, and that it compensates the growth of the volume of $H$. Then the space $\schw_\tau(H)$ (with definition similar to Definition~\ref{def_S}) of smooth functions decreasing $\tau$-rapidly at infinity on $H$, is a Fréchet nuclear space.  
\end{lem}

We know by Proposition~\ref{sous-polynomial} and Proposition~\ref{compense_volume}, that $\sigma$ is sub-polynomial and that it compensates the growth of the volume on $G$. Proposition~\ref{Snuclear} follows from Lemma~\ref{schweitzer}. \end{proof}

\begin{proof}[Proof of Theorem~\ref{Stensor}]
We follow the lines of the proof from the classical result in $\R^n$ \cite[Theorem 51.6]{treves67}. 

First we observe that, thanks to the easy estimate  
\begin{equation} \label{sigma12}
\sigma(g) \, \le \, C \sigma_1(g_1) \sigma_2(g_2), \qquad g=(g_1,g_2) \in G,  \end{equation}
the following inclusion is immediate     
\begin{equation*} \label{inclusionS1xS2}
\SGi\otimes\SGii \subset \SG. 
\end{equation*}

Second, we show that the algebraic tensor product $\SGi\otimes\SGii$ is dense in $\SG$. Indeed, the tensor product $C^\infty_0(G_1)\otimes C^\infty_0(G_2)$ is sequentially dense in $C^\infty(G_1\times G_2)$ \cite[Theorem 39.2]{treves67}, which is itself dense in $\SG$ by Theorem~\ref{DSdense}.  We have then  \[C^\infty_0(G_1)\otimes C^\infty_0(G_2) \subset \SGi\otimes\SGii \subset \SG, \]
with $C^\infty_0(G_1)\otimes C^\infty_0(G_2)$ dense in $\SG$. The density of $\SGi\otimes\SGii$ follows. 

To complete the proof of Theorem~\ref{Stensor}, we have to show that $\SG$ induces the projective topology $\pi$ on $\SGi\otimes\SGii$. We start by showing that it induces a weaker topology. 

To be able to compute, we need to fix some basis $\mathcal{B}_1$, resp. $\mathcal{B}_2$, of the Lie algebra $\g_1$, resp. $\g_2$, and we take $\mathcal{B}= \mathcal{B}_1 \cup \mathcal{B}_2$ as basis of $\g$.  
For $\alpha=(\alpha_1,\alpha_2)\in \N^{m_1}\times \N^{m_2}$, we have 
\[ X_\mathcal{B}^\alpha (\phi \otimes \psi) = X_{\mathcal{B}_1}^{\alpha_1} \phi \otimes X_{\mathcal{B}_2}^{\alpha_2} \psi, \qquad \phi\in \SGi, \, \psi \in \SGii. \]
Then by~(\ref{sigma12}), we have for all $k\in \N$
\[ \| \sigma^k X^\alpha_\mathcal{B} (\phi\otimes\psi) \|_{L^\infty(G)} \le C  \| \sigma_1^k X_{\mathcal{B}_1}^{\alpha_1} \phi \|_{L^\infty(G_1)}
\cdot \| \sigma_2^k X_{\mathcal{B}_2}^{\alpha_2} \psi \|_{L^\infty(G_2)}. \]
It follows that the bilinear mapping $\phi\times\psi\mapsto  \phi\otimes\psi$ from $\SGi\times\SGii$ to $\SGi\otimes\SGii$ equipped with the $\SG$ topology, is separately continuous and hence continuous (since the $\schw_{\sigma_i}(G_i)$ are Fréchet). 

Now by definition, the projective topology $\pi$ on $\SGi\otimes\SGii$ is the strongest topology for which the continuity of the bilinear mapping $\phi\times\psi\mapsto  \phi\otimes\psi$ 
holds. This implies that the topology induced by $\SG$ on $\SGi\otimes\SGii$ is weaker than $\pi$.   

What remains to show now, is 
that $\SG$ induces on $\SGi\otimes\SGii$ a topology   stronger than $\pi$. 

Because the spaces $\schw_{\sigma_i}(G_i)$ ($i=1,2$) are nuclear (Proposition~\ref{Snuclear}), the projective topology $\pi$  coincides with the $\epsilon$-topology on $\SGi\otimes\SGii$.
So to prove that the $\pi$-topology is stronger than the topology inherited from $\SG$, it is enough to show that the identity mapping $I$ from $\SGi\otimes\SGii$ equipped with topology from $\SG$, into $\SGi\otimes\SGii$  equip\-ped with the $\epsilon$-topology, is continuous.  

Since $\SG$ is a Fréchet space, it is metrizable, and so is  $\SGi\otimes\SGii$ equipped with the topology induced by $\SG$. It follows that it is sufficient to show that $I$ is sequentially continuous to deduce that it is continuous. 

On the other hand, if a sequence $\{f_j\}_{j\in \N}$ goes to zero in $\SG$, then it converges uniformly to zero on equicontinuous subset of $\schw'_\sigma(G)$ (by definition of equicontinuity). Let $A'\subset\schw'_{\sigma_1}(G_1)$,  $B'\subset\schw'_{\sigma_2}(G_2)$, be equicontinuous subsets in $\schw'_{\sigma_1}(G_1)$ and $\schw'_{\sigma_1}(G_1)$ respectively. The subset $A'\otimes B'$ 
is equicontinuous in $\schw'_\sigma(G)$. It follows that if $\{f_j\}_{j\in \N} \subset \SGi\otimes\SGii$ goes to zero for the $\SG$-topology, then it converges to zero uniformly on the $A'\otimes B'$. This is valid for every choice of equicontinuous subsets $A'\subset\schw'_{\sigma_1}(G_1)$ and $B'\subset\schw'_{\sigma_2}(G_2)$, so $\{f_j\}_{j\in \N}$ converges to zero for the $\epsilon$-topology. This proves that  the topology induced by $\SG$ is stronger than the $\epsilon$-topology, \emph{i.e.} the $\pi$-topology. Theorem~\ref{Stensor} follows. \end{proof}

As a consequence of Theorem~\ref{Stensor}, we obtain the following generalization of Schwartz Kernel Theorem.

\begin{cor}[Kernel Theorem] \label{skt}
We have the following topological isomorphism
\[L\big(\SGi, \, \schw'_{\sigma_2}(G_2)\big) \cong \schw_\sigma'(G) , \]
where the space of tempered distributions carry the strong dual topology, and the space of continuous linear maps carries topology of uniform convergence on bounded subset. 
\end{cor}

\begin{proof}
Since $\SGi$ and $\SGii$ are nuclear Fréchet spaces, their strong duals are nuclear and we have (for 
the topologies of Corollary~\ref{skt})
\[ L\big(\SGi, \, \schw'_{\sigma_2}(G_2)\big) 
\cong \schw'_{\sigma_1}(G_1) \widehat{\otimes} \schw'_{\sigma_2}(G_2) 
\cong \big( \schw_{\sigma_1}(G_1) \widehat{\otimes} \schw_{\sigma_2}(G_2) \big)', \] 
where $\schw'_{\sigma_1}(G_1) \widehat{\otimes}\schw'_{\sigma_2}(G_2)$ 
denotes the projective completion (which coincides with the $\epsilon$-completion) of the tensor product $\schw'_{\sigma_1}(G_1) {\otimes}\schw'_{\sigma_2}(G_2) $, 
and similarly for $\schw_{\sigma_1}(G_1) \widehat{\otimes} \schw_{\sigma_2}(G_2)$. We refer to \cite[Propositions 50.5 -- 50.7]{treves67} for the proofs of the isomorphisms.  We conclude by Theorem~\ref{Stensor}. \end{proof}

\begin{rem} One can show that the  isomorphism in Corollary~\ref{skt} is defined by \[ \begin{array}[t]{ccc}
L\big(\SGi, \, \schw'_{\sigma_2}(G_2)\big) &\to& \schw_\sigma'(G) \\
\mathcal{K} &\mapsto& K 
\end{array}  \]
with  \[ <K,\phi\otimes \psi>:= <\mathcal{K}\phi, \psi>,\qquad  \phi\in \SGi, \,  \psi \in \SGii.  \]
The distribution $K$ is called the \emph{kernel} of the linear map $\mathcal{K}$. 
\end{rem}


\subsection*{Acknowledgements}
The author wishes to express her gratitude to Andrea Altomani, Carine Molitor-Braun and Fulvio Ricci for several stimulating conversations and suggestions during the preparation of the paper. She whishes  also to thank Michael Cowling for his most helpful 
comments on the weight function $\sigma$.


\bibliographystyle{alpha}

\end{document}